\documentclass[a4paper,10pt]{article}

\usepackage{amssymb,amsmath,graphics}

\newenvironment{dwd}{\par\noindent{\bf Proof.}}{\par\rightline{$\blacksquare$}}

\newtheorem{theo}{Theorem}
\newtheorem{prop}{Proposition}  
\newtheorem{coro}{Corollary}

\newtheorem{rema}{Remark}
\newtheorem{lema}{Lemma}
\newtheorem{defi}{Definition}
\def\be#1\ee{\begin{equation}#1\end{equation}}
\newcommand{\ba}{\begin{eqnarray} }
\newcommand{\ea}{\end{eqnarray} }
\def\bt#1\et{\begin{theo}#1\end{theo}}
\def\bl#1\el{\begin{lema}#1\end{lema}}
\def\bp#1\ep{\begin{prop}#1\end{prop}}
\def\bd#1\ed{\begin{defi}#1\end{defi}}
\global\long\def\sbr#1{\left[ #1\right] }
 \global\long\def\cbr#1{\left\{  #1\right\}  }
 \global\long\def\rbr#1{\left(#1\right)}
 
\def\ccA{{\cal A}}

\def\ccE{{\cal E}}

\def\va{\varepsilon}
\def\ra{\rightarrow}

\def\E{\mathbf{E}}
\def\Cov{\mathbf{Cov}}
\def\P{\mathbf{P}}

\def\R{{\mathbb R}}
\def\Z{{\mathbb Z}}

\def\ls{\leqslant}
\def\gs{\geqslant}

\begin{document}

\title{\bf Sudakov minoration for products of radial-type log concave measures}
\author{{Witold Bednorz}
\footnote{Research partially supported by  NCN Grant UMO-2019/35/B/ST1/04292}
\footnote{Institute of Mathematics, University of Warsaw, Banacha 2, 02-097 Warszawa, Poland}}

\maketitle
\begin{abstract}
The first step to study lower bounds for a stochastic process is to prove a  special property - Sudakov minoration.
The property means that if a certain number of points from the index set are well separated then we can provide an optimal type lower bound for the mean value of the supremum of the process. Together with the generic chaining argument the property can be used to fully characterize the mean value of the supremum of   the stochastic process.
In this article we prove the property for canonical processes based on radial-type log concave measures.
\end{abstract}

\section{Introduction}

Consider a random vector $X$, which takes values in $\R^d$. Let $T$ be a subset of $\R^d$, we can
define canonical process as
\[
(X_t)_{t\in T},\;\;\mbox{where}\;\;X_t=\sum^d_{i=1}t_i X_i=\langle t,X\rangle.
\]
One of the basic questions in the analysis of stochastic processes is to characterize 
$S_X(T)=\E \sup_{t\in T}X_t$, where $(X_t)_{t\in T}$ is a family of 
random variables. It is well known that usually $S_X(T)<\infty$ is equivalent to $\P(\sup_{t\in T}|X_t|<\infty)=1$, i.e that paths are a.s. bounded.
In order to avoid measurablility questions formally
$S_X(T)$  should be defined $\sup_{F\subset T} \E\sup_{t\in F}X_t$, where the supremum is over all finite subsets $F$ of $T$. 
The main tool invented to study the quantity is the generic chaining. Although, the idea works neatly whenever the upper bound is concerned, the lower bound is a much more subtle matter. What one can use towards establishing the lower bound for $\E\sup_{t\in T}X_t$ is the growth condition which can be later used in the partition scheme \cite{TG3} or construction of a special family of functionals \cite{VH}.
The core of the approach is the Sudakov minoration.  Basically,  the minoration means that we can answer the question about 
understanding $\E\sup_{t\in T}X_t$ in the simplest setting, where points in $T$ are well separated.  Before we state the condition, let us make a simple remark.
\begin{rema}\label{s1-rema1}
Let $p\gs 1$.  Suppose that $|T|\ls e^p$, and $\|X_t\|_p\ls A$ for each $t\in T$. Then $\E|\sup_{t\in T}X_t|\ls eA$.
\end{rema}
\begin{dwd}
Indeed,  we have
\[
\E|\sup_{t\in T}X_t|\ls \E\sbr{\sum_{t\in T} |X_t|^p}^{1/p}\ls  \sbr{ \E \sum_{t\in T} |X_t|^p}^{1/p}\ls eA.
\]
\end{dwd}
The Sudakov minoration means that something opposite happen.  Namely,  we require that $|T|\gs e^p$, $p\gs 1$
and for all $s,t\in T$, $s\neq t$, $\|X_t-X_s\|_p\gs A$.  These assumptions should imply
\begin{align}\label{mino-1} 
& \E\sup_{t\in T}X_t\gs K^{-1}A,
\end{align}
where $K$ is an absolute constant. The problem has a long history, which we outline now.
\smallskip

\noindent The property was first proved in 1969 for $X$-Gaussian, i.e. $X=G=(g_i)^d_{i=1}$, $g_i$ are independent standard normal variables.
The result is prescribed to Sudakov \cite{Sud1,Sud2}, although some other researchers could be aware of the fact. In 1986 Pajor and Tomczak-Jaegermann  \cite{P-TJ} observed that the property for Guassian $X$ can be established in the dual way. Then, Talagrand \cite{TG1} proved in 1990, that the minoration works for $X=\va=(\va_i)^d_{i=1}$, $\va_i$ are independent Rademachers. He invented the result, when studying properties of infinitely divisible processes. Later, in 1994 Talagrand \cite{TG2} realized that the minoration holds also for exponential-type distributions. The result was improved by Latala \cite{Lat1} in 1997, so that it holds for $X$ which has independent entries of log-concave tails. The problem to go beyond the class of independent entries occurred to  be unexpectedly hard.  First result towards this direction was due to Latala \cite{Lat2},  2014 and concerned $X$ of density $\exp(-U(\|x\|_p))$, where $U$ is an increasing convex function and $\|x\|_p=\rbr{\sum^d_{i=1}|x_i|^p}^{1/p}$. Then,  in 2015 Latala and Tkocz  \cite{L-T} proved that the minoration in the case of independent entries necessarily  requires that entries are $\alpha$-regular, which basically means that they behaves like variables of log-concave distribution.  Finally,  it has to be mentioned that some new ideas appeared how a proof of the minoration could be established.  For example the dimension reduction - the work of Mendelson,  Milman and Paouris \cite{MPM}.   In this paper there covered some simple cases,  though it is known that the program cannot succeed in full generality.   
\smallskip

\noindent  The aim of the paper is to prove that the minoration works for a certain class of log-concave distributions, that extends mentioned examples.  More precisely,  we show that the property holds for $X$ of density 
\begin{align}\label{mino-2}
\mu_X(dx)=\exp(-\sum^M_{k=1}U_k( \|x_k\|^{p_k}_{p_k} )  )dx,
\end{align} 
where $U_k$ are increasing convex functions, $p_k> 1$ and $x_k=(x_i)_{i\in J_k}$, $J_k$ are pairwise disjoint and covers $[d]=\{1,2,\ldots,d\}$.
We use the notation $n_k=|J_k|$ and $q_k$ for $p_k^{\ast}$, i.e. $1/p_k+1/q_k=1$.
It will be clear from our proof that we can say something about the minoration also if $p_k\gs 1$, however 
if we want the constant $ K $ not to depend on the problem setting, we need some  cutoff level for $ p_k $ from $ 1 $.
\smallskip

\noindent
The distribution $\mu_X$ has some properties that are important for the proof of minoration. 
Namely:
\begin{enumerate}
\item isotropic position;
\item one unconditionality;
\item log-concavity;
\item structural condition.
\end{enumerate}
Let us comment below on them:
\smallskip

\noindent  1) Isotropic position. We say that $X$ is in isotropic position, if
\begin{align}\label{s1-cond0}
\E X=0,\;\;\mbox{i.e.} \E X_i=0,\;\;i\in [d],\;\;\Cov X=\mathrm{Id}_d,\;\;\mbox{i.e.}\;\; \Cov(X_iX_j)=\delta_{i,j}. 
\end{align}
Note that by an affine transformation we can impose the property, keeping the structural assumption \eqref{mino-2}. 
Note also that it is a very natural modification which makes the Eucledean distance
important, i.e. $\|X_t-X_s\|_2=d_2(s,t)=\|t-s\|_2$. 
\smallskip

\noindent 2) We say that $X$ is one-unconditional if
\begin{align}\label{s1-cond1}
& (X_1,X_2,\ldots,X_d)\stackrel{d}{=}(\va_1|X_1|,\va_2|X_2|,\ldots,\va_d|X_d|),
\end{align}
where $|X_i|$, $i\in [d]$ and independent random signs $\va_i$, $i\in [d]$ that
are independent  of all $|X_i|$, $i\in [d]$.  Observe that the assumption
makes it possible to consider $(X_t)_{t\in T}$ in the following form
\[
X_t=\sum^d_{i=1} t_iX_i\stackrel{d}{=}\sum^d_{i=1} t_i\va_i|X_i|.
\] 
3) Log-concavity. We assume that $X$ has the distribution 
\begin{align}\label{s1-cond2}
& \mu(dx)=\exp(-U(x))dx,
\end{align} 
where $U$ is a convex function.  Due to the result of \cite{L-T}  the assumption is almost necessary.
More formally,   what we need is the comparison of moments.  For log-concave and $1$ unconditional $X$,  we have
that  for all $0<p<q$ and $t\in \R^d$  (see \cite{N-O} or  \cite{K-B} for the proof)   
\begin{align}\label{s1-ineq1}
& \|X_t\|_q \ls \frac{\Gamma(q+1)^{\frac{1}{q}}}{\Gamma(p+1)^{\frac{1}{p}}}\|X_t\|_p\ls \frac{q}{p}\|X_t\|_p.
\end{align}
The necessary condition for the Sudakov minoration formulated in \cite{L-T} states that $X$ must be $\alpha$-regular,  which means that there exists $\alpha\gs 1$ such that,  for all $2\ls p<q$ 
\begin{align}\label{s1-ineq2}
& \|X_t-X_s\|_{q}\ls \alpha\frac{q}{p}\|X_t-X_s\|_p.
\end{align}
Clearly,  for log-concave and  $1$ unconditional  $X$ the inequality \eqref{s1-ineq2} works with $\alpha=1$.
\smallskip

\noindent
4) Structural assumption.  As we have explained in \eqref{mino-2} we assume that the distribution of $X$ has 
a certain structure.  More precisely, we require that
\begin{align} \label{s1-cond3}
& U(x)=\sum^M_{k=1} U_k(\|x_k\|^{p_k}_{p_k}),\;\;p_k>1.
\end{align}
Therefore we can treat $X$ as $(X_k)^M_{k=1}$,  where vectors $X_k$ with values in $\R^{n_k}$ are independent. 
The fact, which we use later is that each $X_k$ has the same distribution as $R_kV_K$, where $R_k$ and $V_k$ are independent and
$R_k$ is distributed on $\R_{+}$ with the density $x^{n-1}\exp(-U_k(x^{p_k}))|\partial B^{n_k}_{p_k}|$ and $V_k$ is distributed on $\partial B^{n_k}_{p_k}$
with respect to the (probabilistic) cone measure $\nu_k$  - for details see e.g. \cite{Woj}. The most important property
of the cone measure is that for any integrable $f:\R^{n_k}\ra\R$
\[
\int_{\R^{n_k}} f(x) dx=|\partial B^{n_k}_{p_k}| \int_{\R_{+}} \int_{\partial B^{n_k}_{p_k}} f(r\theta) r^{n_k-1} \nu_k(d\theta) dr.
\]
Note that here $B^{n_k}_{p_k}=\{x\in \R^{n_k}:\;\|x\|_{p_k}\ls 1\}$ and
$\partial B^{n_k}_{p_k}=\{x\in \R^{n_k}:\;\|x\|_{p_k}=1\}$.

\section{Results}

Following the previous section, we can formulate the main result of this paper
\begin{theo}\label{diety-sudmin}
Suppose that $X\in \R^d$ is a random vector whose density is of the form
$\mu_X=\mu_1\otimes \mu_2 \otimes \ldots \otimes \mu_M$ and $\mu_k$ is a log concave distribution on
$\R^{n_k}$ given by 
\[
\mu_k(dx_k)=\exp(-U_k(\|x_k\|^{p_k}_{p_k}))dx_k,\;\;\mbox{where}\;\;x_k=(x_{i})_{i\in J_k},
\]
where $J_k$ are disjoint sets such that $\bigcup^M_{k=1}J_k=[d]$ and $p_k\gs 1+\va$ for some $\va>0$.
Then,  the Sudakov minoration holds,  i.e.  for every set $T\subset \R^d$ such that $\|X_t-X_s\|_p\gs A$,  $s\neq t$,  $s,t\in T$
for some $p\gs 1$ the following inequality holds
\[
\E\sup_{t\in T}X_t \gs K^{-1}A,
\]
 where $K$ depends on $\va$ only.  
\end{theo}
The proof is quite complicated that is why it is good to give a sketch of our approach.
\begin{enumerate}
\item Simplifications.  We show that Theorem \ref{diety-sudmin} has to be established for $A=p$,  $X$ in the isotropic position and sets  $T$ such that for each $t\in T$,  $t=(t_i)^d_{i=1}$,  $t_i\in \{0,k_i\}$ for some positive $k_i$,  $i\in [d]$ and $\|X_t-X_s\|_p\gs p$,  $s\neq t$, $s,t\in T$.  Moreover,  for each $t\in T$,   $\sum^d_{i=1}|t_i|$ and $\sum^d_{i=1}1_{t_i=k_i}$ are much smaller than $p$.  In other words,  minoration should be proved for cube-like sets with thin supports.
\item Moments.  The next step is a careful analysis of the condition $\|X_t-X_s\|_p\gs p$.   In particular,  we are going to use the structure assumption $X=(X_1,X_2,\ldots,X_M)$,  where entries  $X_k\in \R^{n_k}$ are independent.  The main trick here is to define a random vector $Y=(Y_1,Y_2,\ldots,Y_M)$ such that $Y$ has all coordinates  independent,  i.e.
not only $Y_k\in \R^{n_k}$ are independent but also coordinates of each $Y_k$ are independent.  Due to our structure 
assumption it will be possible to represent $X_k=R_k V_k$,  $Y_k=\tilde{R}_k V_k $,  where $R_k,\tilde{R}_{k}\in \R_{+}$
are some $\alpha$ regular variables and $V_k$ are distributed with respect to cone measures on $\partial B^{n_k}_{p_k}$ respectively.
\item Split.  We observe that we can split each point $t\in T$ into $t^{\ast}$ - small part and $t^{\dagger}$ - large part.
More precisely,  we decide whether $t^{\ast}_k$ is $t_k$ or $0$ and respectively $t^{\dagger}_k$ is $0$ or $t_k$ depending on
how large some norm of $t_k$ is. 
Moreover,  we show that either there is a subset $S\subset T$,  $|S|\gs e^{p/2}$  such that
$\| X_{t^{\ast}}-X_{s^{\ast}}\|_p\gtrsim p$  for any $s\neq t$ $s,t\in S$ or we  can find
 $S\subset T$,  $|S|\gs e^{p/2}$ such that $\|X_{t^{\dagger}}-X_{s^{\dagger}}\|_{p}\gtrsim p$ for any $s\neq t$,  $s,t$ in $S$.
\item Small part.  If the 'small part' case holds then we show that one may forget about variables $R_k,\tilde{R}_k$ and
in this way deduce the minoration from such a result for the random vector $Y$.
\item Large part.  If the 'large part' case holds,  then we prove that variables $V_k$ are not important,  more precisely 
we prove,  following the approach from the 'simplifications' step,  that not only $t_i\in \{0,k_i\}$ but we have such a property  for  
$t^{\dagger}_k=(t^{\dagger}_i)_{i\in J_k}$,  namely $t^{\dagger}_k\in \{x_k,0\}$ for some $x_k\in \R^{n_k}$.  Consequently,  we 
deduce the minoration from such a result for the random vector  $(X_{x_k})^M_{k=1}$.
 \end{enumerate}
It is a standard argument that having the minoration one can prove the comparison between
$\E \sup_{t\in T}X_t$ and $\gamma_X(T)$.  However,  this approach requires some two additional 
conditions.  First,  we need that there exists $\va>0$ such that
\begin{align}\label{s1-ineq3}
& \|X_t\|_{2^{n+1}}\gs (1+\va)\|X_t\|_{2^{n}},\;\;n\gs 0.
\end{align}
The property works if there is a cutoff level for all $p_k$ below $\infty$. 
Moreover,  there must hold a certain form of measure concentration.
What suffices is that there exist constants $K,L\gs 1$
such that for any $p\gs 2$ and set $T\subset \R^d$
\begin{align}\label{s1-ineq4}
& \left\| \rbr{\sup_{t\in T}X_t -K \E \sup_{t\in T} X_t}_{+} \right\|_p\ls L\sup_{s,t\in T}\|X_t-X_s\|_p.
\end{align}
The result is known only in few cases.  Fortunately,  in our setting it can be derived from the infimum convolution \cite{La-Wo}.  Note that the problem can be easily reduced  to the one where $T$ is a gauge of some norm in $\R^d$.  Moreover,
the main concentration inequality -  $CI(\beta)$ from \cite{La-Wo} holds for radial log-concave densities and it can be  tensorized,  namely  if $CI(\beta)$ holds for measures $\mu_1,\mu_2,\ldots,\mu_M$ then $CI(\beta')$ holds for
$\mu_1\otimes \mu_2 \otimes \ldots \otimes \mu_M$ with some dependence between $\beta$ and $\beta'$,  i.e.   $\beta'$ is some multiplication of $\beta$.   
\begin{theo}\label{diety-character}
Under the assumption of Theorem \ref{diety-sudmin} and assuming additionally that $p_k\ls p_{\infty}<\infty$
then
\[
\E\sup_{t\in T}X_t \simeq \gamma_X(T).
\]
\end{theo}
The result is quite standard and goes through the idea of growth condition - \cite{TG3}. 
\smallskip

\noindent
Before we start the proof of Theorem \eqref{diety-sudmin} we need a preliminary result which explains that the Sudakov minoration has to be established only for sets $T$ that have cube-like structure.

\section{Cube-like sets}

First step concerns some basic simplifications of the problem. We first note that the minoration has to be proved only for certain cube-like sets $T$.
We use symbols $\lesssim,\gtrsim,\simeq$ whenever we compare quantities up to a numerical constant comparable to $1$.
 \smallskip
 
 \noindent
We have to start from  the Bobkov Nazarov \cite{B-N} inequality or rather from its 
basic consequence: 
\begin{align}\label{s2-ineq1}
& \|X_t-X_s\|_p \lesssim \|\ccE_t-\ccE_s\|_p, 
\end{align}
where $(\ccE_i)^d_{i=1}$ are independent symmetric standard exponential variables.
The Bobkov Nazarov inequality concerns one-unconditional and log concave distributions in the isotropic position.
However,  the inequality \eqref{s2-ineq1} may hold in a bit more general setting.  For example,  it is also true
if entries of an isotropic $X$ are independent and $\alpha$-regular. 
 That is why we simply refer to \eqref{s2-ineq1}
as the basic requirement.  In the  proof of the simplification,  the condition is used together with the Kwapien-Gluskin result \cite{K-G} 
\begin{align}\label{s2-ineq2}
& \|\ccE_t\|_p \sim p\|t\|_{\infty}+\sqrt{p}\|t\|_2.
\end{align}
On the other hand,  by the result of Hitczenko \cite{Hit},  we have also the lower bound
\begin{align}\label{s2-ineq3}
& \|X_t\|_p \gtrsim \|\va_t\|_p\sim \sum_{1\ls i\ls p}|t^{\ast}_i|+\sqrt{p}\rbr{\sum^d_{i>p}|t_i^{\ast}|^2}^{1/2},
\end{align}
where $(t^{\ast}_i)^d_{i=1}$ is the rearrangement of $t$ such that $|t^{\ast}_1|\gs |t^{\ast}_2|\gs \ldots \gs |t^{\ast}_d|$.
The last tool we need is that Sudakov minoration works for the random vector $(\va_i)^d_{i=1}$.  As we have already  mentioned the result was first  proved in \cite{TG1}.  We are ready to formulate the main simplification result
\begin{prop}\label{s2-prop1}
Suppose that $X$ is in the isotropic position and
Fix $p\gs 2$. Suppose that:
\begin{enumerate}
\item  $X$ is one unconditional;
\item  $X$ is log concave;  
\item  $X$ satisfies \eqref{s2-ineq1}. 
\end{enumerate}
Then to show the Sudakov minoration it suffices to prove 
for all sets $T$ such that:
\begin{enumerate}
\item  $\exp(p)\ls |T|\ls 1+\exp(p)$ and $0\in T$;
\item  for each $i\in [d]=\{1,2,\ldots,d\}$
\[
t_i\in \{0,k_i\},\;\;\mbox{where}\;k_i\gs \rho,
\]
where $\rho\ls e^{-1}$ and $\rho/\log \frac{1}{\rho}= 4C\delta$;
\item  for each $t\in T$
\begin{align}\label{s2-ineq4}
\sum_{i\in I(t)} k_i \ls 2C\delta p,\;\;\mbox{where}\;\;I(t)=\{i\in [d]:\;t_i=k_i\};
\end{align}
\item  for all $s,t\in T$, $s\neq t$
\[
\|X_t-X_s\|_p\gs p,
\] 
\end{enumerate}
where $\delta$ is suitably small  and $C\gs 1$ a universal constant,
that the following inequality holds
\[
\E\sup_{t\in T} X_t\gs K^{-1}p,
\]
where $K$ is a universal constant.
\end{prop}
\begin{dwd}
The proof is based on a number of straightforward steps. 
\smallskip

\noindent 
{\bf Step 1} Obviously it suffices to show that the Sudakov minoration works for $A=p$.
Let $m_p(t,s)=\|\sum^n_{i=1} (t_i-s_i)\va_i\|_p$.
We may assume that $p\gs 1$ is suitably large.  Moreover, we may consider $T$ such that $0\in T$, $|T|\gs  \exp(\frac{p}{2})+\exp(\frac{3p}{4})$ and
$m_p(t,0)\ls \delta p$ for all $t\in T$,   where $\delta\ls 1$ can be suitably small.
\smallskip

\noindent
By the Talagrand's result \cite{TG2} (see also \cite{TG3} for the modern exposition)  
if $N(T,m_p,u)\gs \exp(\frac{p}{4})-1$ then $\E \sup_{s,t\in T} \sum^n_{i=1}(t_i-s_i)\va_i\gs L^{-1} u$,  for
a universal $L$.  Therefore either the result holds with $K=2^{-1} L^{-1}  \delta$ or  
$\E \sup_{s,t\in T} \sum^n_{i=1}(t_i-s_i)\va_i< 2^{-1}L^{-1} \delta p$ and then
\[
N(T,m_p,\frac{1}{2}\delta p)\ls \exp(\frac{p}{4})-1\ls \frac{\exp(\frac{p}{2})}{1+\exp(\frac{1}{4}p)}.
\]
It implies that there exists $t_0\in T$ such that 
\[
\{t\in T:\; m_{p}(t,t_0)\ls \delta p\}|\gs |T|\frac{1+\exp(\frac{1}{4}p)}{\exp(\frac{p}{2})}\gs \exp(\frac{p}{2})+\exp(\frac{3p}{4}).
\]
Therefore,  we may consider set $T'=\{t-t_0: m_{p}(t,t_0)\ls \delta p\}$,  which satisfies
all the requirements we have promised. 
\smallskip

\noindent
{\bf Step 2} 
Let $\delta\ls e^{-1}$ and  $\rho/\log(1/\rho)=4C\delta$.
We may assume that $0\in T$,  $|T|\gs 1+\exp(\frac{p}{4})$ and additionally  
\[
t_i\in (k_i-\rho,k_i+\rho)\cup (-\rho,\rho )\;\;\mbox{for all}\;t\in T\;\;\mbox{and}\;\; i\in [d],
\] 
where $k_i$ are given numbers such that $k_i\gs \rho$ where $\rho\ls e^{-1}$ and satisfies 
\[
\rho\log \frac{1}{\rho}=4 \delta\ls e^{-1},\;\;\mbox{where}\; C \gs 1
\] 
is a universal constant which we choose later. 
\smallskip

\noindent
Indeed.  consider measure  $\mu=\otimes^d_{i=1}\mu_i$, where $\mu_i(dx)=\frac{1}{2}e^{-|x|}dx$ for
all $i\in [d]$.  For any $x\in \R^d$,  $x=(x_i)^d_{i=1}$ and $t\in T$,  $t=(t_i)^d_{i=1}$
\[
T_x=\{t\in T:\;t_i\in (x_i-\rho,x_i+\rho)\cup (-\rho,\rho),\;i\in [d]\}
\]
and
\[
A_t=\{x\in \R^d:\;t_i\in (x_i-\rho,x_i+\rho)\cup (-\rho,\rho),\;i\in [d]\}.
\]
Now there are two possibilities either
\begin{align}\label{s2-ineq4.5}
\mu_i(\{x_i:\;t_i\in (x_i-\rho,x_i+\rho)\cup (-\rho,\rho) \})\gs \rho e^{-|t_i|-\rho}
\end{align}
or
\begin{align}\label{s2-ineq4.6}
\mu_i(\{x_i:\;t_i\in (x_i-\rho ,x_i+\rho)\cup (-\rho,\rho) \})=1,\;\;\mbox{if}\;|t_i|<\rho.
\end{align}
Applying \eqref{s2-ineq3} we get for some $C\gs 1$
\begin{align}\label{s2-ineq5} 
& m_{p}(t,0)=\|\sum^d_{i=1}t_i\va_i\|_{p}\gs C^{-1}\rbr{\sum_{i\ls p}|t^{\ast}_i|+\sqrt{p}\rbr{\sum_{i>p}|t_i^{\ast}}^2}^{\frac{1}{2}},
\end{align}
where $|t^{\ast}_i|$ is the non-decreasing rearrangement of $|t_i|$.  Again we  deduce from \eqref{s2-ineq3} that
\begin{align}\label{s2-ineq6} 
& m_{p}(t,0)\ls C\rbr{\sum_{i\ls p}|t^{\ast}_i|+\sqrt{p}\rbr{\sum_{i>p}|t_i^{\ast}}^2)}^{\frac{1}{2}},
\end{align}
for some $C_2\gs 1$.  Therefore,  using \eqref{s2-ineq5} and $\rho/\log \frac{1}{\rho}=4C\delta$ we obtain 
\begin{align}\label{s2-ineq7} 
|\{i\in \{1,...,d\}:\;|t_i|\gs \rho\}|\ls \frac{p}{4\log \frac{1}{\rho}}\ls \frac{p}{4}.
\end{align}
Again using \eqref{s2-ineq5} this shows 
\begin{align}\label{s2-ineq8}
\sum^d_{i=1} |t_i| 1_{|t_i|\gs \rho}\ls C\delta p.
\end{align}
Consequently,  by (\ref{s2-ineq4.5},\ref{s2-ineq4.6})
\[
\mu(A_t)\gs \rho^{\frac{p}{4\log \frac{1}{\rho}}}\exp(-2C \delta p)  \gs \exp(-\frac{p}{2}).
\]
However,  using that $|T|\gs \exp(\frac{3p}{4})$ we infer
\[
\int \sum_{t\in T}1_{A_t}(x)\mu(dx)\gs |T|\exp(-\frac{p}{2} )\gs \exp(\frac{p}{4}).
\]
Therefore we get that there exists
at lest one point $k\in \R^d$ such that
\[
|T_k|\gs \exp(\frac{p}{4}).
\]  
It is obvious that $|k_i|$ may be chosen in a way that $|k_i|\gs \rho$.
Combining \eqref{s2-ineq8}  with $|k_i|\gs \rho$ and $t_i\in (k_i-\rho,k_i+\rho)$,  we obtain
\[
\sum^d_{i=1} (|k_i|-\rho)\vee \rho 1_{|t_i|\gs \rho}\ls C\delta p.
\]
Clearly $(|k_i|-\rho)\vee \rho\gs \frac{1}{2}|k_i|$ and therefore
\[
\frac{1}{2} \sum^d_{i=1}|k_i|1_{|t_i|\gs \rho}\ls C \delta p,
\]
which implies \eqref{s2-ineq4}.  Clearly by the symmetry of each $X_i$ we may only consider positive $k_i\gs \rho$. 
\smallskip

\noindent
{\bf Step 3} It suffices to consider set $T$ which additionally satisfies $t_i\in \{0,k_i\}$ where $k_i\gs \rho$.  Moreover, $\exp(p/4)\ls T\ls 1+\exp(p/4)$,  $0\in T$ 
and 
\begin{align}\label{s2-ineq9}
& \|X_t-X_s\|_p\gs \frac{p}{2},\;\;\mbox{for all}\;s,t\in T,s\neq t.
\end{align}
Consider the following function
\[
\varphi_i(t_i)=\left\{\begin{array}{lll}
0 & \mbox{if} & |t_i|<\rho \\
k_i & \mbox{if} & |t_i|\gs \rho 
\end{array}\right.  .
\]
Let $\varphi(t)=(\varphi_i(t_i))^d_{i=1}$. We show that $\| X_{\varphi(t)}-X_{\varphi(s)}\|_p\gs \frac{p}{2}$. 
It requires  some upper bound on $\|X_{t-\varphi(t)}\|_p$.  Consider any $s\in T$ then using \eqref{s2-ineq2}
\[
\| X_s\|_p\ls  C'(p\|s\|_{\infty}+\sqrt{p}\|s\|_2).
\]
Note that for $s=t-\varphi(t)$ we get by the contraction principle $m_p(s,0)\ls p\rho+m_p(t,0)$ and hence using \eqref{s2-ineq6}
\begin{align*}
&\|X_s\|_p\ls  C'(p\|s\|_{\infty}+\sqrt{p}\|s\|_2)\ls 2C'\rho p+CC' m_p(s,0) \\
&\ls C'(2+C)\rho p+CC' m_p(t,0))\ls C'((2+C)\rho+C\delta)p\ls \frac{p}{4},
\end{align*}
for suitably small $\delta$ and hence also suitably small $\rho$.  Therefore
\[
\|X_{\varphi(t)}-X_{\varphi(s)}\|_{p}\gs \|X_{t}-X_{s}\|_{p}-\|X_{t}-X_{\varphi(t)}\|_{p}-\|X_{s}-X_{\varphi(s)}\|_{p}\gs \frac{p}{2}.
\]
Suppose we can prove the main result for the constructed set $T$,  which satisfies $|T|\gs \exp(p/4)$.  Formally,  we select a subset $S\subset T$
such that $0\in S$ and $\exp(p/4)\ls |S|\ls 1+\exp(p/4)$ .  The Sudakov minoration for cube-like sets gives
\[
\E \sup_{t\in S}X_{\varphi(t)}\gs K^{-1}p,
\]  
for some universal $K$.  Recall that 
\[
\|X_{t-\varphi(t)}\|_p\ls C'((2+C)\rho+C\delta)p
\]
and therefore by Remark \ref{s1-rema1}  and $\exp(p/4)\ls |S|\ls 1+\exp(p/4)$,  $0\in S$ we get
\[
\E\sup_{t\in S} X_{t-\varphi(t)}\ls eC'\rbr{(2+C)\rho+C\delta}p.
\] 
Thus
\[
\E \sup_{t\in S} X_t=\E \sup_{t\in S} X_{\varphi(t)}+X_{t-\varphi(t)}\gs \E \sup_{t\in S}X_{\varphi(t)}- \E\sup_{t\in S} X_{t-\varphi(t)}\gs \frac{1}{2}K^{-1}p,
\]
for suitably small $\delta$,  i.e.  $2eC' \rbr{(2+C)\rho+C\delta}\ls K^{-1}$.  Obviously,  the set $S$ in the required simplification
in this step.
\smallskip

\noindent
{\bf Step 4} The final step is to replace $p$ by $4p$ so that $\exp(p)\ls |T|\ls 1+\exp(p)$ and $0\in T$.
Then obviously by\eqref{s1-ineq2} (with $\alpha=1$)
\[
4 \|X_t-X_s\|_{p}\gs \|X_t-X_s\|_{4p}\gs 2p.
\]
Thus we may to redefine $t_i$ as $2 t_i$.  Consequently,  the theorem holds with slightly
rearranged constants,  namely we set $\delta'=2 \delta$ instead of $\delta$ and  $\rho'$ (instead of $\rho$) that satisfies $\rho'/\log(1/\rho')= 4C\delta'$.  Obviously $\rho'\ls e^{-1}$ if
$\delta$ was suitably small. 
\end{dwd}
The above proof is a slightly rearranged version of the argument presented in \cite{Lat1} - we have stated the proof here for the sake of completeness.  Note also that without much effort the argument works for $\alpha$ regular $X$,  i.e..
when the inequality \eqref{s1-ineq2} is satisfied for all $s,t\in \R^d$. 
\begin{rema}\label{s2-rem0}
There are suitably small $\delta'$ and $\delta''$ such that for any set $T$ that satisfies properties from Preposition \ref{s2-prop1} 
then for all $t\in T$ the following holds 
\begin{align}\label{s2-ineq10}
|I(t)|\ls \delta' p,\;\; \|t\|_1=\sum_{i\in I(t)}k_i \ls \delta'' p,
\end{align}
\end{rema}
\begin{dwd}
 Indeed by \eqref{s2-ineq4} we get $\|t\|_1\ls 2C\delta p$.
On the other hand,  $\rho |I(t)|\ls \|t\|_1\ls 2C\delta p$.  Since  $4C\delta/\rho= 1/\log(1/\rho)$ this implies
\[
|I(t)| \ls \frac{1}{2} \frac{1}{\log \frac{1}{\rho}}p.
\] 
We set $\delta'=2C\delta$ and $\delta''=1/(2\log(1/\rho))$.
\end{dwd}
There is another property which we can add to our list of conditions that $T$ has to satisfy. Namely, 
it suffices to prove minoration only if  $\|X_t-X_s\|_p\simeq p$ for all $t\in T$.
\begin{prop}\label{s2-prop2}
Suppose that random vector $X$ satisfies \eqref{s1-ineq2}, then it suffices to prove the minoration only for sets $T$ such that
$p\ls \|X_t-X_s\|_p\ls 2\alpha p$.
\end{prop}
\begin{dwd}
The argument is rather standard and can be found e.g.  in \cite{Le-Ta}.  Basically, either there is at least $e^{p/2}$ points that
are within the distance $2\alpha p$ from some point in $T$ or one can find at least  $e^{p/2}$ points that are $2\alpha p$ separated, i.e. $\|X_t-X_s\|_p\gs 2\alpha p$
for all $s,t$ $s\neq t$ in the set.  However,  then $\|X_t-X_s\|_{p/2}\gs \frac{2\alpha}{2\alpha}p=p$, then $\|X_t-X_s\|_{p/2}\gs p$. 
We continue with this set instead of $T$.  It is easy to understand that in this way we have to find a subset $T'$ of $T$ that counts at least
$e^{p/2^m}$ elements,  where $2^m\ls p$,  such that 
\[
p\ls \|X_t-X_s\|_{p/2^m}\ls 2\alpha p\;\; s,t\in T'.
\]
In this case obviously, for $ \tilde{t}=t/2^m$ and $\tilde{T}=\{\tilde{t}:\; t\in T'\}$ we get $\|X_{\tilde{t}}-X_{\tilde{s}}\|_{p/2^m}\simeq p/2^m$. Then, by our assumption that the minoration works for $\tilde{T}$ and $p/2^m$
\[
\frac{p}{K2^m}\ls \E\sup_{\tilde{t}\in \tilde{T}}X_{\tilde{t}}=\E\sup_{t\in T'}X_{t/2^m}\ls \frac{1}{2^m}\E\sup_{t\in T}X_t.
\]
Otherwise,  if we reach $m$ such that $2^{m}\ls p< 2^{m+1}$,  we obtain that there are at least two points $s,t$, $s\neq t$ such that
 $\|X_t-X_s\|_2\gs  p $. This immediately implies the minoration,  since
\[
\E\sup_{t\in T}X_t\gs \frac{1}{2}\E |X_t-X_s|\gs \frac{1}{2\sqrt{2}} \|X_t-X_s\|_2\gtrsim p.
\]
Obviously,  by homogeneity
\[
\E\sup_{t\in T}X_t\gs \E\sup_{t\in \bar{T}}X_{t/2^m} =\frac{1}{2^m} \E\sup_{t\in T'}X_t \gs K^{-1}p,
\]
which ends the proof.
\end{dwd}
There is one more useful remark.  As explained in previous works on the problem - see Lemma 2.6 in \cite{Lat2} or section 4 in  \cite{MPM} - the case when $p\gs d$ is 'easy'.  More precisely,  
\begin{rema}\label{s2-rem1}
If Sudakov minoration holds for $p=d$,  then it holds also for $p\gs d$.
\end{rema}

\section{How to compute moments}

The second step concerns basic facts on moments of $X_t$,  $t\in T$.
Recall that we work with a random vector  $X=(X_1,X_2,\ldots,X_d)$ which is $1$ unconditional,  isotropic and log-concave.    In particular it implies that  for any $t\in \R^d$
\[
X_t=\langle X, t \rangle =\sum^d_{i=1}t_i X_i \stackrel{d}{=}\sum^d_{i=1} t_i \va_i|X_i|,
\]
where $\va=(\va_i)^d_{i=1}$ is a vector of independent Rademacher  variables which is independent of $X$.
We start from a series of
general facts that are known in this case. 
We also discuss moments of  $X$ with independent,  symmetric and $\alpha$-regular entries.
It will be discussed later that due to our basic simplification - Proposition \ref{s2-prop1} - we have to compute 
$\|X_t\|_p$  only when  $d\ls p$,  which is  a bit simpler than the general case. 
We start from the characterization  proved in \cite{Lat2}. 
\begin{theo}\label{s3-theo1}
Suppose that $X$ has $1$-unconditional and log-concave distribution $\mu_X(dx)=\exp(-U(x))dx$.
We assume also that $X$ is in the isotropic position.  Then for any $p\gs d$
\begin{align}\label{s3-ineq1}
& \|X_t\|_p\simeq \sup \cbr{\sum^d_{i=1}|t_i| a_i:\;\; U(a)-U(0) \ls p }.
\end{align}
\end{theo}
In \cite{Lat2},  there is also an alternative formulation of this result
\begin{theo}\label{s3-theo2}
Under the same assumptions as in Theorem \ref{s3-theo1} we have for any $p\gs d$
\begin{align}\label{s3-ineq2}
& \|X_t\|_p\simeq \sup \cbr{\sum^d_{i=1}|t_i| a_i:\;\;\P\rbr{ \bigcap^d_{i=1}\{|X_i|\gs a_i\}}\gs e^{-p} }.
\end{align}
 \end{theo}
Note that a similar characterization works when $X$ has independent and $\alpha$-regular entries.
There is also  a version of the above fact formulated in terms of moments.   Namely,
\begin{theo}\label{s3-theo3}
Under the same assumptions as in Theorem \ref{s3-theo1} we have for any $p\gs d$
\begin{align}\label{s3-ineq4}
& \|X_t\|_p\simeq \sup\cbr{\sum^d_{i=1}|t_i| \|X_i\|_{a_i}: \;\;\sum^d_{i=1}a_i\ls p }.
\end{align}
\end{theo}
\begin{dwd}
We prove the result for the sake of completeness. 
For simplicity we consider only log-concave, one unconditional,  isotropic case.  The $\alpha$-regular case can be proved in the similar way.
Consider $a_i\gs 0$, $i\in [d]$ such that $\sum^d_{i=1}a_i\ls p$.  Since $X$ is isotropic $a_i<2$ are unimportant,  that is why we use $a_i\vee 2$.
It is well-known that moments and quantiles  are comparable,  namely for some constant $C\gs 1$  
\[
\P\rbr{ |t_i X_i|\gs C^{-1}\|t_i X_i\|_{a_i\vee 2}}\gs C^{-1}e^{-a_i\vee 2}
\]
and hence 
\begin{align}\label{s3-ineq55}
\prod^d_{i=1}\P\rbr{|t_iX_i|\gs C^{-1}\|t_iX_i\|_{a_i\vee 2}} \gs 
C^{-d}e^{-\sum^d_{i=1}a_i\vee 2}.
\end{align}
Consequently,
\begin{align*}
& \|X_t\|_p\gs C^{-2}\rbr{\sum^d_{i=1}\|t_iX_i\|_{a_i\vee 2}}e^{-\sum^d_{i=1}\frac{a_i\vee 2}{p}}\\
& \gs C^{-2}e^{-3}\sum^d_{i=1}\|t_iX_i\|_{a_i\vee 2},
\end{align*}
where we have used that
\[
\sum^d_{i =1}a_i\vee 2\ls 2p+\sum^d_{i=1}a_i\ls 3p.
\]
We turn to prove the converse inequality. 
Let $|X|=(|X_i|)^d_{i=1}$ and $|X|_t=\sum^d_{i=1}|t_i||X_t|$. Since $p\gs d$, $\|X\|_p\ls \||X|_t\|_p\ls 2\|X_t\|_p$.
By the homogeneity of 
the problem,  we can assume that $\||X|_t\|_p=p$ (possibly changing point $t$ by a constant).
We have to prove that there exists $a_i$ such that $\sum^d_{i=1}a_i\ls p$
whereas $\sum^d_{i =1}\|t_iX_i\|_{a_i}\gs C^{-1}p$. It is clear that if suffices to prove 
the result for $p$ suitably large in particular for $p\gs 2$. 
Let $\gamma$ be a constant which we determine later. We define $r_i$ as 
\begin{enumerate}
\item  $r_i=2$ if $\|t_iX_i\|_2\ls  2\gamma $
\item  $r_i=p$ if $\|t_i X_i\|_p\gs p\gamma$ 
\item  or otherwise $r_i=\inf\{r\in [2,p]:\;\|t_iX_i\|_r=r\gamma\}$ 
\end{enumerate}
We first observe that if $r_{i_0}\gs p$ then there is nothing to prove since choosing
 $a_{i_0}=r_{i_0}$ and other $a_i$ equal $0$ we fulfill the requirement 
$r_{i_0}=\sum^d_{i =1}a_i\ls p$ whereas 
\[
\sum^d_{i =1}\|t_iX_i\|_{a_i\vee 2}
\gs \|t_iX_i\|_{a_{i_0}}\gs \gamma a_{i_0}\gs\gamma p.
\]
Therefore, we may assume that $r_i< p$ for all $i\in [d]$.
Now suppose that we can find a subset $J\subset [d]$ with the property $\sum_{i\in J}r_i\gs 3p$
and $J$ does not contain a smaller subset with the property.  Therefore necessarily
$\sum_{i\in J}r_i\ls 4p$ since $r_i\ls p$ for any $i\in [d]$.  Obviously,
\[
\sum^d_{i =1}\|t_iX_i\|_{r_i}\gs \sum_{i\in J}\|t_iX_i\|_{r_i}\gs 
\gamma \sum_{i\in J}(r_i -2)\gs  \gamma p
\] 
and hence we can set $a_i=r_i/4$, which implies that $\sum^d_{i=1} a_i\ls p$ and
\[
4\sum^d_{i =1}\|t_i X_i\|_{a_i\vee 2}\gs \sum^d_{i=1}\|t_i X_i\|_{r_i}\gs \gamma p.
\]
Therefore,  to complete the proof we just need to prove that $\sum^d_{i=1}r_i\gs 4p$. 
By the log-concavity for $r_i\gs 2$, $u\gs 1$
\[
\P\rbr{t_i|X_i|\gs  e u \|t_i X_i\|_{r_i}}\ls e^{-r_i u} 
\]
and therefore 
\[
\P\rbr{\frac{t_i |X_i|}{\gamma  e}\gs r_i u }\ls e^{-r_i u}=\P(|\ccE_i|\gs r_i u),
\]
where we recall that $\ccE_i$ are symmetric standard exponentials.
Consequently,
\begin{align*}
&\||X|_t\|_p\ls \rbr{ e \gamma \sum^d_{i =1}   r_i}+\\
&+\sbr{\E \rbr{\sum^d_{i =1}|t_i| |X_i|1_{|t_iX_i|\gs \gamma  e r_i}}^p }^{\frac{1}{p}}\ls \\
&\ls  ( e\gamma  \sum^d_{i =1} r_i)+ 
\gamma e\rbr{\E |\sum^d_{i =1} |\ccE_i|1_{|\ccE_i|\gs r_i}|^p }^{\frac{1}{p}} \\
&\ls \rbr{ e \gamma \sum^d_{i =1}  r_i}+\gamma e \|Z\|_p,
\end{align*}
where $Z$ is of gamma distribution $\Gamma(d,1)$. 
Clearly $\|Z\|_p\ls (p+d)\ls  2p$ and therefore
\[
p=\||X|_t\|_p\ls \gamma e \sum^d_{i=1} r_i+4\gamma  e p. 
\]
Therefore choosing $\gamma=1/(8e)$
\[
\frac{1}{2}\||X|_t\|_p\ls \frac{1}{8}\sum^d_{i=1}r_i.
\]
It proves that $\sum^d_{i=1}r_i\gs 4p$,  which completes the proof.
\end{dwd}
We  turn to prove some remarks on moments in our case..  We are interested in moments of
$\|X_t\|_p$ for $X$ that satisfies our structural assumption.   In order to explain all the idea we need a lot of random variables.  It is good to collect their definitions in one place in order to easily find the right reference.
\bd
We define the following random variables.
\begin{itemize}\label{s3-defi0}
\item Let $X\in \R^d$ be a random vector in the isotropic position.
\item Let $X=(X_1,X_2,\ldots,X_M)$ where $X_k\in \R^{n_k}$,  $k\in [M]$ are independent and $X_k$ has the density
\[
\exp(-U_k(\|x_k\|^{p_k}_{p_k})),\;\;\mbox{where}\;\; U_k \;\; \mbox{is\;\; convex, \;\; inceasing},\;\; x_k=(x_{i})_{i\in J_k}.
\]
\item Sets $J_k$,  $1\ls k\ls M$ are disjoint and $|J_k|=n_k$. 
\item We denote $X_k=(X_i)_{i\in J_k}$,  we also use the ordering $X_k=(X_{k1},X_{k2},\ldots,X_{kn_k})$.
\item Random variables $R_k\in \R_{+}$,  $V_k\in \R^{n_k}$,  $k\in [M]$,  are such that $X_k=R_kV_k$.
\item Random vector $V_k=(V_{k1},V_{k2},\ldots,V_{kn_k})$ is distributed with respect to the cone measure on $\partial B^{n_k}_{p_k}$.
\item Random variable $R_k$ is distributed on $\R_{+}$ with respect to the density
\begin{align}\label{s3-ineq5}
& g_k(s)=s^{n-1}\exp(-U_k(s^{p_k}))|\partial B^{n_k}_{p_k}| 1_{\R_{+}}(s).
\end{align} 
\item Let $Y=(Y_1,Y_2,\ldots,Y_M)$,    where $Y_k\in \R^{n_k}$,  $k\in [M]$ are independent and $Y_k$ has the density
\begin{align}\label{s3-eq1}
f_k(x_k)=\prod_{i\in J_k } \frac{b_k^{\frac{1}{p_k}}}{2\Gamma\rbr{1+\frac{1}{p_k}}}
e^{-b_k|x_i|^{p_k}}, \;\;\mbox{where}\;\; b_k=\sbr{\frac{\Gamma(\frac{3}{p_k})}{\Gamma(\frac{1}{p_k})}}^{\frac{p_k}{2}}.
\end{align}
\item We denote $Y_k=(Y_i)_{i\in J_k}$,  we also use the ordering $Y_k=(Y_{k1},Y_{k2},\ldots,Y_{kn_k})$.
\item Random variables $\tilde{R}_k\in \R_{+},V_k \in \R^{n_k}$,  are such that $Y_k=\tilde{R}_k V_k$,  moreover 
$V_k$ is already defined and distributed like the cone measure on $\partial B^{n_k}_{p_k}$.
\item Random variable $\tilde{R}_k$ is distributed 
on $\R_{+}$ with respect to the density
\begin{align}\label{s3-ineq5.5}
\tilde{g}_k(s)= \frac{b_k^{\frac{n_k}{p_k}}}{\frac{1}{p_k}\Gamma(\frac{n_k}{p_k})} s^{n_k-1}\exp(-b_k s^{p_k}) 1_{\R_{+}}(s),\;\;b_k=\sbr{\frac{\Gamma(\frac{3}{p_k})}{\Gamma(\frac{1}{p_k})}}^{\frac{p_k}{2}}.
\end{align}
\end{itemize}
\ed
There are some basic consequences of these definitions.  Note that two of them,  i.e.  formulas 
(\ref{s3-ineq5},  \ref{s3-ineq5.5})
we have mentioned above. 
\begin{rema}\label{s3-rema0} 
There is a list of basic properties of variables described in Definition \ref{s3-defi0}.
\begin{itemize}
\item Variables $R_k$,  $k\in [M]$ and $V_k$,  $k\in [M]$ are all independent of each other.
\item Variables $\tilde{R}_k$,  $k\in [M]$ and $V_k$,  $k\in [M]$ are all independent of each other.
\item All variables $Y_i$,  $i\in J_k$,  $k\in [M]$ are independent and isotropic,  in fact $Y_i$,  $i\in J_k$
has the density
\[
 \frac{b_k^{\frac{1}{p_k}}}{2\Gamma\rbr{1+\frac{1}{p_k}}}
e^{-b_k|x|^{p_k}},\;\;x\in \R.
\] 
\item Variables $V_k$, $k\in [M]$ are $1$ unconditional. 
\item We have $|\partial B_{p_k}^{n_k}|=p_k\rbr{2\Gamma\rbr{1+1/p_k}}^{n_k}/ \Gamma( n_k / p_k )$.
\item Clearly $X_{t_k}=R_k\langle V_k, t_k \rangle$ and hence $\|X_{t_k}\|_r=\|R_k\|_r\|\langle V_k,t_k\rangle \|_{r}$,  $r>0$.  In particular 
\begin{align}\label{s3-ineq5.6}
1=\E X_{ki}^2=\E R_k^2 \E  V_{ki}^2, \;\;i=1,2,\ldots,n_k. 
\end{align}
\item In the same way $Y_{t_k}=\tilde{R}_k \langle V_k,t_k \rangle$ and hence $\|Y_{t_k}\|_r=\|\tilde{R}_k\|_r\|\langle V_k,t_k\rangle \|_{r}$,  $r>0$.  In particular 
\begin{align}\label{s3-ineq5.7}
1=\E Y_{ki}^2= \E \tilde{R}_k^2 \E  V_{ki}^2, \;\;i=1,2,\ldots,n_k.
\end{align}
\end{itemize}
\end{rema}
Let $I_k(t)=\{i\in J_k:\; |t_i|>0\}$.  Using Theorem \ref{s3-theo1} we get for $r\gs |I_k(t)|$, where
\begin{align}
\label{s3-ineq5.8} \|\langle Y_k,t_k \rangle\|_r=\|\tilde{R}_k\|_r\|\langle V_k,t_k \rangle\|_r 
\simeq r^{1/p_k}\|t_k\|_{q_k}.
\end{align}
On the other hand,
\[
\|\tilde{R}_k\|_r  =b^{-\frac{1}{p_k}}_k\frac{\Gamma\rbr{\frac{n_k+r}{p_k}}^{\frac{1}{r}}}{\Gamma\rbr{\frac{n_k}{p_k}}^{\frac{1}{r}}}\simeq \rbr{n_k+r}^{\frac{1}{p_k}}.
\]
Consequently,
\begin{prop}\label{s3-prop1}
The following holds
\[
 \|\langle V_k,t_k \rangle\|_r \simeq\left\{ \begin{array}{lll}
 \|t_k\|_{q_k} \frac{r^{\frac{1}{p_k}}}{n_k^{\frac{1}{p_k}}} & \mbox{if} &
 |I_k(t)|\ls r\ls n_k \\
 \|t_k\|_{q_k}   &  \mbox{if} & r>n_k \end{array}
\right.
\]
\end{prop}
It remains to explain the role of $\|R_k\|_r$.
For all $k\in [M]$ we define 
\begin{align}\label{s3-ineq5.9}
S_k=R_k^{p_k}, \;\; \tilde{U}_k\;\;\mbox{satisfies}\;\;  \exp(-\tilde{U}_k(x))1_{\R_{+}}(x)=\exp(-U_k(x))\frac{1}{p_k}|\partial B^{n_k}_{p_k}|1_{\R_{+}}(x).
\end{align}
Clearly  $S_k$ has the density
\[
h_k(x)=x^{\frac{n_k}{p_k}-1}\exp(-U_k(x))\frac{1}{p_k}|\partial B^{n_k}_{p_k}|1_{\R_{+}}(x)=x^{\frac{n_k}{p_k}-1}\exp(-\tilde{U}_k(x))1_{\R_{+}}(x).
\]
Note that $\|R_k\|_{r}=\|S_k\|^{\frac{1}{p_k}}_{r/p_k}$.
We use the result of Ball \cite{K-B}-Lemma 4.
\begin{lema}\label{s3-lema1}
For any $0<p<q$, the following holds
\begin{align}\label{s3-ineq6}
& e^{-\tilde{U}(0)q} \Gamma(p+1)^{q+1}
\sbr{\int_{\R_{+}}  x^{q}e^{-\tilde{U}(x)} dx }^{p+1}
\ls e^{-\tilde{U}(0)p} \Gamma(q+1)^{p+1} 
\sbr{\int_{\R_{+}}  x^{p}e^{-\tilde{U}(x)} dx }^{q+1},
\end{align}
where $\tilde{U}(x)$ is a  convex,  increasing function on $\R_{+}$.
\end{lema}
Note that we use the above result for $\tilde{U}_k$ defined in \eqref{s3-ineq5.9}.
Consequently,  we get
\begin{coro}\label{s3-coro1}
For any $r\gs 1$, we have 
\[
\|R_k\|_r=\|S_k\|^{\frac{1}{p_k}}_{\frac{r}{p_k}}\ls e^{\frac{1}{n_k}\tilde{U}_k(0)} \frac{\Gamma\rbr{\frac{n_k+r}{p_k}}^{\frac{1}{r}}}{\Gamma\rbr{\frac{n_k}{p_k}}^{\frac{1}{r}+\frac{1}{n_k}}}.
\]
\end{coro}
We have to bound $\exp(-\tilde{U}_k(0)/n_k)$. 
Fist using the isotropic position  of $Y_k$ we get
\[
1= \E V_{ki}^2 \E \tilde{R}_k^2 =\frac{1}{b_k^{\frac{2}{p_k}}}\frac{\Gamma\rbr{\frac{n_k+2}{p_k}}}{\Gamma\rbr{\frac{n_k}{p_k}}} \E V_{ki}^2,
\]
hence
\[
\E V_{ki}^2=\frac{b_k^{\frac{2}{p_k}}\Gamma\rbr{\frac{n_k}{p_k}}}{\Gamma\rbr{\frac{n_k+2}{p_k}}}. 
\]
Therefore, by Lemma \ref{s3-lema1}   we get
\[
\frac{1}{\E (V_{ki})^2}=\E R_k^2 \ls e^{\frac{2}{n_k}\tilde{U}_k(0)}\frac{\Gamma\rbr{\frac{n_k+2}{p_k}}}{\Gamma\rbr{\frac{n_k}{p_k}}^{1+\frac{2}{n_k}}}.
\]
It proves the following inequality
\begin{align}\label{s3-ineq7}
& e^{\frac{1}{n_k}\tilde{U}_k(0)}\gs b_k^{-\frac{1}{p_k}}\Gamma\rbr{\frac{n_k}{p_k}}^{\frac{1}{n_k}}.
\end{align}
In order to prove the upper bound we need a Hensley type inequality.
First observe that
\[
 \P(R_k\ls  t)=\int^t_0 x^{n_k-1}e^{-\tilde{U}_k(x^{p_k})}dx\ls \int^t_0 x^{n_k-1}e^{-\tilde{U}_k(0)} dx
 =\frac{1}{n_k}t^{n_k}e^{-\tilde{U}_k(0)}.
\]
Let $F_k(t)=\frac{1}{n_k}t^{n_k}e^{-\tilde{U}_k(0)}$ for $0\ls t\ls t_{\ast}$,
where $t_{\ast}=\rbr{n_k e^{\tilde{U}_k(0)}}^{1/n_k}$.
We can calculate
\begin{align*}
& \frac{1}{\E (V_k)^2_i}=\E R_k^2= 2\int^{\infty}_0 t\P(R_k>t) dt\\
& \gs  2\int^{t_{\ast}}_0 t(1-F_k(t))dt\gs 2\int^{t_{\ast}}_0 t(1-t^{n_k})dt\\
&=t_{\ast}^2-2\frac{1}{n_k+2}t_{\ast}^{n_k+2}e^{-\tilde{U}_k(0)}\frac{1}{n_k} \\
& =t_{\ast}^2-\frac{2}{n_k+2}t_{\ast}^2=\frac{n_k}{n_k+2}t_{\ast}^2.
\end{align*}
Therefore
\[
\rbr{n_k e^{\tilde{U}_k(0)}}^{\frac{2}{n_k}}\ls \frac{n_k+2}{n_k} \frac{\Gamma\rbr{\frac{n_k+2}{p_k}}}{b_k^{\frac{n_k}{p_k}}\Gamma\rbr{\frac{n_k}{p_k}}}.
\]
Which means
\begin{align}\label{s3-ineq8}
& e^{\frac{1}{n_k}\tilde{U}_k(0)} \ls  \frac{1}{n_k^{\frac{1}{n_k}}}\frac{(n_k+2)^{\frac{1}{2}}}{n_k^{\frac{1}{2}}}
\frac{\Gamma\rbr{\frac{n_k+2}{p_k}}^{\frac{1}{2}}}{b_k^{\frac{1}{p_k}}\Gamma\rbr{\frac{n_k}{p_k}}^{\frac{1}{2}}}.
\end{align}
It finally proves
\[
\frac{\Gamma\rbr{\frac{n_k}{p_k}}^{\frac{1}{n_k}} }{b_k^{\frac{1}{p_k}}}\ls
e^{\frac{1}{n_k}\tilde{U}_k(0)} \ls \frac{1}{n_k^{\frac{1}{n_k}}}\frac{(n_k+2)^{\frac{1}{2}}}{n_k^{\frac{1}{2}}}
\frac{\Gamma\rbr{\frac{n_k+2}{p_k}}^{\frac{1}{2}}}{b_k^{\frac{1}{p_k}}\Gamma\rbr{\frac{n_k}{p_k}}^{\frac{1}{2}}}
\]
and hence
\[
e^{\frac{1}{n_k}\tilde{U}_k(0)} \simeq \frac{1}{\|V_{ki}\|_2}\simeq n_k^{\frac{1}{p_k}}. 
\]
Together with Corollary \ref{s3-coro1} it shows that
\begin{align}\label{s3-ineq9}
& \|R_k\|_r\lesssim n_k^{\frac{1}{p_k}}\frac{\Gamma\rbr{\frac{n_k+r}{p_k}}^{\frac{1}{r}}}{\Gamma\rbr{\frac{n_k}{p_k}}^{\frac{1}{n_k}+\frac{1}{r}}}.
\end{align}
Note that
\[
\frac{\Gamma\rbr{\frac{n_k+r}{p_k}}^{\frac{1}{r}}}{\Gamma\rbr{\frac{n_k}{p_k}}^{\frac{1}{n_k}+\frac{1}{r}}}\simeq \rbr{1+\frac{r}{n_k}}^{\frac{n_k+r}{p_k r}}.
\]
That is why
\begin{align}\label{s3-ineq10}
& \|R_k\|_r\lesssim (n_k+r)^{1/p_k}.
\end{align}
In particular $\|R_k\|_r \simeq n_k^{1/p_k}$ for $1\ls r\ls n_k$.
We introduce variables $P_k=R_k/\|R_k\|_2$,  $Q_k=\tilde{R}_k/\|\tilde{R}_k\|_2$, $W_k=\|R_k\|_2 V_k$.
This helps since now vectors $(\va_k P_k)^M_{k=1}$, $(\va_k Q_k)^M_{k=1}$ 
are in the isotropic  position,  where $\va_k$, $k=1,2,\ldots,M$
are independent Rademacher variables,  independent of all $P_k$ and $Q_k$.
Obviously,  also $W_k$,  $k=1,2,\ldots, M$ are isotropic,  independent and independent of all $P_k$.
Thus in particular $\|P_k\|_r,\|Q_k\|_r \simeq 1$, for $1\ls r\ls n_k$.
Once again due to Lemma \ref{s3-lema1}
\begin{align*}
& e^{-\tilde{U}_k(0)(\frac{n_k+q}{p_k}-1)} \Gamma\rbr{\frac{n_k+p}{p_k}}^{\frac{n_k+q}{p_k}}
\sbr{\int_{\R_{+}}  x^{\frac{n_k+q}{p_k}-1}e^{-\tilde{U}_k(x)} dx }^{\frac{n_k+p}{p_k}}\\
& \ls e^{-\tilde{U}_k(0)(\frac{n_k+p}{p_k}-1)} \Gamma\rbr{\frac{n_k+q}{p_k}}^{\frac{n_k+p}{p_k}} 
\sbr{\int_{\R_{+}}  x^{\frac{n_k+p}{p_k}-1}e^{-\tilde{U}_k(x)} dx }^{\frac{n_k+q}{p_k}}.
\end{align*}
Therefore,
\[
\|R_k\|_q^{\frac{q(n_k+p)}{p_k}}\ls e^{\tilde{U}_k(0)\frac{(q-p)}{p_k}}\frac{\Gamma\rbr{\frac{n_k+q}{p_k}}^{\frac{n_k+p}{p_k}}}{\Gamma\rbr{\frac{n_k+p}{p_k}}^{\frac{n_k+q}{p_k}}}\|R_k\|_p^{\frac{p(n_k+q)}{p_k}}. 
\]
which is
\begin{align}
& \|R_k\|_q \ls e^{\tilde{U}_k(0) \frac{(q-p)}{q(n_k+p)}}\frac{\Gamma\rbr{\frac{n_k+q}{p_k}}^{\frac{1}{q}}}{\Gamma\rbr{\frac{n_k+p}{p_k}}^{\frac{n_k+q}{q(n_k+p)}}}
\|R_k\|_p^{\frac{p(n_k+q)}{q(n_k+p)}}. 
\end{align}
Note that
\[
 \frac{\Gamma\rbr{\frac{n_k+q}{p_k}}^{\frac{1}{q}}}{\Gamma\rbr{\frac{n_k+p}{p_k}}^{\frac{n_k+q}{q(n_k+p)}}}
\simeq \frac{(n_k+q)^{\frac{n_k+q}{p_kq}}}{(n_k+p)^{\frac{n_k+q}{p_kq}}}  \simeq \rbr{\frac{n_k+q}{n_k+p}}^{\frac{n_k+q}{qp_k}}.
\]
Moreover, for $p\gs 1$
\[
e^{\tilde{U}_k(0) \frac{(q-p)}{q(n_k+p)}}\lesssim (n_k^{\frac{1}{p_k}})^{\frac{(q-p)n_k}{q(n_k+p)}}\lesssim \|R_k\|_p^{\frac{(q-p)n_k}{q(n_k+p)}}
\]
and since 
\[
\frac{(q-p)n_k}{q(n_k+p)}+\frac{p(n_k+q)}{q(n_k+p)}=1
\]
it proves that
\begin{align}\label{s3-ineq12}
& \|R_k\|_q \lesssim \rbr{\frac{n_k+q}{n_k+p}}^{\frac{1}{p_k}}\|R_k\|_p,
\end{align}
which means that $R_k$ is $\alpha$ regular.  In particular $\va_k R_k$, $k=1,2,\ldots,M$ are independent symmetric
$\alpha$-regular variables.  It also implies that $\tilde{R}_k$ have the fastest growth among all the possible distributions of  $R_k$.
\smallskip

\noindent We end this section with the characterization of $\|X_t-X_s\|_p$. 
Recall our notation $P_k=R_k /\|R_k\|_2$ and $W_k=\|R_k\|_2 V_k$.
We show 
\begin{prop}\label{s3-prop2}
Under the assumptions as  in Theorem \ref{s3-theo1}  and assuming  that all the properties mentioned in Preposition \ref{s2-prop1} are satisfied,  the following result holds
\begin{align}\label{s3-ineq13}
& \|X_t-X_s\|_{p}\simeq \sup\cbr{ \sum^M_{k=1} \|P_k\|_{r_k}\| \langle W_k,t_k-s_k
\rangle\|_{r_k}1_{r_k\gs |I_k(t)\triangle I_k(s)||}:\;\;\sum^M_{k=1}r_k=p } \\
\label{s3-ineq13.5} &  \|Y_t-Y_s\|_p \simeq  \sup\cbr{ \sum^M_{k=1} \|Q_k\|_{r_k}\| \langle W_k,t_k-s_k
\rangle\|_{r_k}1_{r_k\gs |I_k(t)\triangle I_k(s)||}:\;\;\sum^M_{k=1}r_k=p }.
\end{align}
Moreover,  if $r_k\gs |I_k(t)\triangle I_k(s)||$ then 
\[
\|Q_k\|_{r_k}\| \langle W_k,t_k-s_k\rangle \|_{r_k}=\|Y_{t_k}\|_{r_k}\sim r_k^{1/p_k}\|t_k\|_{q_k} \lesssim \|P_k\|_{r_k}| \langle W_k,t_k-s_k\rangle\|_{r_k}=\|X_{t_k}\|_{r_k}.
\]
\end{prop}
\begin{dwd}
We use that $X_t-X_s=\sum^M_{k=1} \langle X_k, t_k-s_k \rangle$ and $X_k$ are independent,  log-concave vectors.  Observe that 
\[
X_{t_k-s_k}= \langle X_k, t_k-s_k \rangle=P_k \langle W_k, t_k-s_k\rangle.
\] 
Let us denote also $I_k(t)=J_k\cap I(t)$.
By Proposition \ref{s2-prop1} - see \eqref{s2-ineq10} - the number of $k\in M$ for which $t_k-s_k$ is non-zero is much smaller than $p$.   We use Theorem \ref{s3-theo3}  treating $X_t-X_s$ as $\sum^M_{k=1}X_{t_k-s_k}$,  it yields  
\[
\|X_t-X_s\|_{p}\simeq \sup\cbr{\sum^M_{k=1} \|X_{t_k}\|_{r_k}: \;\;\sum^M_{k=1}r_k=p }.
\]
Since $\|X_{t_k-s_k}\|_{r_k}=\|P_k\|_{r_k}\| \langle W_k,  t_k-s_k\rangle\|_{r_k}$ we get
\begin{align}\label{s4-ineq0}
& \|X_t-X_s\|_{p}\simeq \sup\cbr{\sum^M_{k=1} \|P_k\|_{r_k}\| \langle W_k,t_k-s_k
\rangle\|_{r_k}:\;\;\sum^M_{k=1}r_k=p }.
\end{align}
Moreover,
\[
\| P_k \|_r \simeq 1,\;\;\mbox{for}\;\;r\ls n_k,\;\;\|P_k\|_r \lesssim \frac{r_k^{\frac{1}{p_k}}}{n_k^{\frac{1}{p_k}}},\;\;r>n_k.
\]
and 
\[
\| \langle W_k,t_k-s_k\rangle\|_{r}\simeq 
\left\{ \begin{array}{lll}  r^{1/p_k}\|t_k-s_k\|_{q_k}  & \mbox{for} & |I_k(t)\triangle I_k(s)|\ls r\ls n_k \\
n_k^{1/p_k}\|t_k-s_k\|_{q_k}  & \mbox{for} & r>n_k         \end{array}\right.  .
\]
We stress that we should not care about  $r_k<|I_k(t)\triangle I_k(s)|$
since by \eqref{s2-ineq10},  $\sum^M_{k=1}|I_k(t)\triangle I_k(s)|$ is much smaller than $p$. Thus $r_k<|I_k(t)\triangle I_k(s)|$
can be ignored when estimating $\|X_t-X_s\|_p$.  Thus finally 
\[
\|X_t-X_s\|_{p}\simeq \cbr{\sum^M_{k=1} \|P_k\|_{r_k}\| \langle W_k,t_k-s_k
\rangle\|_{r_k}1_{r_k\gs |I_k(t)\triangle I_k(s)||}:\;\;\sum^M_{k=1}r_k=p }
\]
as required.  The proof for $\|Y_t-Y_s\|_p$ follows the same scheme,  we have to use \eqref{s3-ineq5.8}.
\end{dwd}
We will use the above result many times in the subsequent proofs,  even not mentioning it.  

\section{Positive process}

We are going to show that it suffices to prove that for all $p\gs 1$
\[
\E \sup_{t\in T}\sum^M_{k=1}|X_{t_k}|\gs \frac{1}{K}p.
\]
For simplicity we define $|X|_t=\sum^M_{k=1}|X_{t_k}|$.
We show that
\begin{lema}\label{s4-lema1}
Suppose that $T$ satisfies simplifications from Proposition \ref{s2-prop1}. Then 
$\E \sup_{t\in T}|X|_t\gs \frac{1}{K}p$,  implies that also 
$\E \sup_{t\in T} X_t\gs \frac{1}{4K}p$.
\end{lema}
\begin{dwd}
Before we start,  let introduce notation: writing $\E_I$ for $I\subset [M]$ we mean the integration over vectors $(X_k)_{k\in I}$.
Note in particular that $\E=\E_{I,I^c}$.
\smallskip

\noindent
Let us observe that
\begin{align*}
 &\E\sup_{t\in T} X_t=\E\sup_{t\in T}\sum^M_{k=1}X_{t_k}\\
 &=\E\sup_{t\in T}\sum^M_{k=1}\va_k X_{t_k}\gs \E\sup_{t\in T}\sum^M_{k=1}\va_k |X_{t_k}|,
\end{align*}
where $(\va_k)^M_{k=1}$ are independent random signs. The last inequality is due to Bernoulli comparison
and the inequality
\[
||X_{t_k}|-|X_{s_k}||\ls |X_{t_k}-X_{s_k}|.
\] 
Now the standard argument works,  namely
\begin{align*}
& \E\sup_{t\in T}\sum^M_{k=1}\va_k |X_{t_k}|=\sum_{I\subset [M]}\frac{1}{2^M}\E_{I,I^c} \sup_{t\in T} \sum_{k\in I}|X_{t_k}|-\sum_{k\in I^c}|X_{t_k}|\\
&\gs \sum_{I\subset [M]}\frac{1}{2^M}\E_{I} \sup_{t\in T}\rbr{\sum_{k\in I}|X_{t_k}|-\sum_{k\in I^c}\E_{I^c}|X_{t_k}|}\\
& \gs \sum_{I\subset [M]}\frac{1}{2^M}\E \sum_{k\in I}|X_{t_k}|-\delta" p\gs \frac{1}{2}\E\sup_{t\in T}\sum^M_{k=1}|X_{t_k}|-\delta" p,
\end{align*}
where we have used Jensen's inequality and \eqref{s2-ineq10} which implies
\[
\sum_{k\in I^c}\E |X_{t_k}|\ls \sum_{i\in I(t)} k_i\E|X_i|\ls \sum_{i\in I(t)} k_i \ls \delta" p.
\]
It suffices that $\delta"\ls 1/(4K)$.
\end{dwd}
The above result enables us to split points $t_k$ into small and large part.  More precisely,  let us split $t_k=t^{\ast}_k+t^{\dagger}_k$,  where
\begin{align}
\label{s4-ineq2.1} {\bf small\;\; part:}\;\;  t^{\ast}_k=t_k\;\; \mbox{if} \;\; D^{q_k}\|t_k\|_{q_k}^{q_k}\ls A^{p_k}n_k\;\; \mbox{and} \;\;0\;\; \mbox{otherwise}, \\
\label{s4-ineq2.2} {\bf large\;\; part:}\;\; t^{\dagger}_k=t_k \;\; \mbox{if} \;\; D^{q_k}\|t_k\|_{q_k}^{q_k}>A^{p_k} n_k \;\; \mbox{and}\;\; 0\;\; \mbox{otherwise}.
\end{align}
We will need that $A,D\gs 1$ and $A$ is  such that $A^{p_k}\gs D^{q_k}$ for each $k\in [M]$.
This is the moment when we take advantage of  the cutoff level $p_k\gs 1+\va$, since then we can find  such $A$ that does depend on $\va$ only.  In fact we could do better since our split is not important when $q_k>n_k$. Indeed, then
$n_k^{1/q_k}\simeq 1$,  $\|t_k\|_{q_k}\gs \|t_k\|_{\infty}$, and it is well known that we can require in the proof of minoration
that each $|t_i|$, $i\in I(t)$ is large enough. That means,  for $p_k$ too close to $1$,  the first case simply cannot happen.  
In general,  we do not want to bound $n_k$ but fortunately,  it is possible to reduce the dimension of 
our problem by a simple trick.  Note that $T$ has about $e^p$ points whose supports are thin - much smaller than $p$. Thus,  there
are only $p e^p$  important coordinates on which we can condition our basic vector $X$.   Moreover,
due to the  Prekopa-Leindler theorem,  such a reduced problem is still of log-concave distribution
type.   Since it does not affect computation of $\|X_t-X_s\|_p$ moments $s,t$ we end up in the question, 
where $d\ls pe^p$,  which together with Remark \ref{s2-rem1} show that we should care only for $p\ls d\ls pe^p$. 
The trouble with this concept is that it introduces dependence of $\va$ on p which we do not like in
our approach.
\smallskip

\noindent
We have two possibilities:
\begin{itemize}
\item either there exists $S\subset T$, $|S|\gs e^{p/2}$  such that for any $s,t\in S$, $s\neq t$ and
\begin{align}\label{s4-ineq3}
\sum^M_{k=1}D^{q_k}\|t_k^{\ast}-s^{\ast}_k\|_{q_k}^{q_k}=\sum^M_{k=1} \sum_{i\in I_k(t^{\ast})\triangle I_k(s^{\ast})} D^{q_k}k_i^{q_k}\gs p ;
\end{align}
\item or there is a point $t_0\in T$ and a  subset $S\subset T$, $|S|\gs e^{p/2}$ such that for any $t\in S$ we have
\begin{align}\label{s4-ineq4}
\sum^M_{k=1}D^{q_k}\|t_k^{\ast}-t^{\ast}_{0k}\|_{q_k}^{q_k}=\sum^M_{k=1} \sum_{i\in I_k(t^{\ast})\triangle I_k(t_{0}^{\ast})} D^{q_k}k_i^{q_k}< p.
\end{align}
\end{itemize}
Let $\E_P,\E_W$ denote the integration with respect to $(P_k)^M_{k=1}$ and $(W_k)^M_{k=1}$.
Note that we have
\begin{align}
\nonumber & \E \sup_{t\in T}|X|_t\gs \E \sup_{t\in T}\sum^M_{k=1}|P_k||\langle W_k,t_k \rangle|\\
\label{s4-ineq4.5} & \gs \E_W \sup_{t\in T} \sum^M_{k=1}\E_P |P_k| |\langle W_k,t_k \rangle|\simeq \E \sup_{t\in T}\sum^M_{k=1} |\langle W_k,t_k \rangle|. 
\end{align}
For simplicity we denote
\[
|\tilde{X}|_t=\sum^M_{k=1}|\langle W_k,t_k\rangle|.
 \]
Obviously, $|\tilde{X}|_t\gs |\tilde{X}|_{t^{\ast}}$.  We are going to prove that in the case of small coefficients (cf. \eqref{s4-ineq2.1}) necessarily 
$\E\sup_{t\in T}|\tilde{X}|_{t^{\ast}}\gs K^{-1}p$.

\section{Small coefficients}

Recall that we work in the cube-like setting introduced in Proposition \ref{s2-prop1}.
We are going to prove the Sudakov minoration, in the setting where
there are a lot of well separated points in $T$, in the sense  of \eqref{s4-ineq3}.
Toward this goal, we need the process $(Y_t)_{t\in T}$ - see \eqref{s3-eq1}. We slightly modify the process, namely let
\[
Z_t=\sum^M_{k=1} \sum_{i\in I_k(t)} \va_i \rbr{k_i|Y_i|}\wedge D^{q_k}k_i^{q_k},\;\;|Z|_t=\sum^M_{k=1} \sum_{i\in I_k(t)} \rbr{k_i|Y_i|}\wedge D^{q_k}k_i^{q_k},
\]
where as usual $\va_i$, $i\in [d]$ are independent random signs.
Since now the entrances are independent, in order to prove that $\E\sup_{t\in S} Z_{t^{\ast}}\gs K^{-1}p$ it suffices
to show that  $\|Z_{t^{\ast}}-Z_{s^{\ast}}\|_p\gtrsim p/D$.  We prove that this holds under the condition  \eqref{s4-ineq3}.
\begin{lema}\label{s5-lema1}
Suppose that \eqref{s4-ineq3} holds.  Then $\|Z_{t^{\ast}}-Z_{s^{\ast}}\|_p\gtrsim p/D$.
\end{lema}
\begin{dwd}
Note that for some $\sum^M_{k=1}\sum_{i\in I_k}r_{ki}=p$
\[
\|Z_{t^{\ast}}-Z_{s^{\ast}}\|_p\simeq \sum^M_{k=1}\sum_{i\in I_k(t^{\ast})\triangle I_k(s^{\ast})}\|\rbr{k_i|Y_i|}\wedge D^{q_k}k_i^{q_k}\|_{r_{ki}}.
\]
However,  it is clear that $\|\rbr{k_i|Y_i|}\wedge D^{q_k}k_i^{q_k}\|_{r_{ki}}\simeq \|k_iY_i\|_{r_{ki}}\wedge D^{q_k}k_i^{q_k}$.
Since $\|k_iY_i\|_{r_{ki}}\simeq k_ir_{ki}^{\frac{1}{p_k}}$ we get
\[
\|Z_{t^{\ast}}-Z_{s^{\ast}}\|_p\simeq \sum^M_{k=1}\sum_{i\in I_k(t^{\ast})\triangle I_k(s^{\ast})}\rbr{k_ir_{ki}^{\frac{1}{p_k}}}\wedge D^{q_k}k_i^{q_k}.
\]
Our assumption is that
\[
\sum^M_{k=1}\sum_{i\in I_k(t^{\ast})\triangle I_k(s^{\ast})}D^{q_k} k_i^{q_k}\gs p,
\]
which means that we can select  sets $I_k\subset I_k(t^{\ast})\triangle I_k(s^{\ast})$ such that  $r_{ki}=D^{q_k}k_i^{q_k}$ for each  $i\in I_k$
and $r_{ki}=0$ for $i\in J_k\backslash I_k$,  which satisfy  $\sum^M_{k=1}\sum_{i\in I_k}D^{q_k}k_i^{q_k}\simeq p$.  Then
\[
\|Z_{t^{\ast}}-Z_{s^{\ast}}\|_p\gtrsim \sum^M_{k=1}\sum_{i\in I_k}D^{q_k-1}k_i^{q_k}\gtrsim p/D.
\]
This proves the result.
\end{dwd}
We have proved that $\E\sup_{t\in S}Z_{t^{\ast}}\gs K^{-1}p$.
This implies that
\[
K^{-1} p \ls  \E\sup_{t\in S} \sum^M_{k=1} \rbr{\sum_{i\in I_k(t^{\ast})} k_i\va_i |Y_i|} \wedge D^{q_k}\|t_k\|_{q_k}^{q_k}.
\]
The crucial thing is to establish a similar inequality replacing $\big|\sum_{i\in I_k(t)} k_i \va_i |Y_i  |\big|$ with
$\sum_{i\in I_k(t)}k_i|Y_i|$.  This can be done following the approach we have used in the proof of the positive process lemma - Lemma \ref{s4-lema1}.  Namely,  we have
\begin{lema}\label{s5-lema2}
The following inequality holds
\begin{align}\label{s5-ineq0}
 \E\sup_{t\in S} \sum^M_{k=1} \rbr{\sum_{i\in I_k(t^{\ast})} k_i \va_i |Y_i  | }\wedge D^{q_k}\|t_k\|_{q_k}^{q_k}\gs K^{-1}p.
 \end{align}
\end{lema}
\begin{dwd}
Let us observe that
\begin{align*}
& \E\sup_{t\in S} \sum^M_{k=1}|\sum_{i\in I_k(t^{\ast})} k_iY_i \big| \wedge \|t_k\|_{q_k}^{q_k}\\
&=\frac{1}{2^d} \sum_{I\subset [d]} \E_{I,I^c}\sup_{t\in S} \rbr{\sum^M_{k=1} \big| \sum_{i\in I_k(t^{\ast})\cap I} k_i|Y_i|-\sum_{i\in I_k(t^{\ast})\cap I^c}k_i|Y_i| \big| \wedge D^{q_k}\|t_k\|_{q_k}^{q_k}}\\
&\gs \frac{1}{2^d} \sum_{I\subset [d]} \E_I \sup_{t\in S} \rbr{\sum^M_{k=1} \rbr{\sum_{i\in I_k(t^{\ast})\cap I} k_i|Y_i|} \wedge D^{q_k}\|t_k\|_{q_k}^{q_k}-
\E_{I^c}\sum_{i\in I_k(t^{\ast})\cap I^c}k_i|Y_i|} \\
&\gs \frac{1}{2^d} \sum_{I\subset [d]} \E\sup_{t\in S} \sum^M_{k=1} \rbr{\sum_{i\in I_k(t^{\ast})\cap I} k_i|Y_i|} \wedge D^{k_q}\|t_k\|_{q_k}^{q_k}-\delta'' p\\
& \gs\frac{1}{2}\E\sup_{t\in S} \sum^M_{k=1} \rbr{\sum_{i\in I_k(t^{\ast} } k_i|Y_i|} \wedge D^{q_k}\|t_k\|_{q_k}^{q_k}-\delta'' p,
\end{align*}
where we have used $\E |Y_i|\ls 1$ and $\sum_{i\in I(t)}k_i\ls \delta'' p$ - Proposition \ref{s2-prop1},  see \eqref{s2-ineq10}.  
We have used also that for  positive $a,b,c$,  inequalities $|a-b|\wedge c\gs a\wedge c - b$ and $a\wedge c+b\wedge c\gs (a+b)\wedge c$.
For suitably small $\delta''$ this implies 
\[
(4K)^{-1}p \ls \E \sup_{t\in S} \sum^M_{k=1} \big|\sum_{i\in I_k(t^{\ast})} k_i \va_i |Y_i  | \big| \wedge D^{q_k}\|t_k\|^{q_k}_{q_k}.
\]
\end{dwd}
Let us recall our notation $Q_k=\tilde{R}_k/\|\tilde{R}_k\|_2$,  $W_k=V_k \|\tilde{R}_k\|_2$,  $Q_k\langle W_k, t_k\rangle =Y_{t_k}$.
Consequently,  the above result can be rewritten in the following form
\begin{align}\label{s5-ineq0.5}
(4K)^{-1}p \ls \E \sup_{t\in S} \sum^M_{k=1} \rbr{Q_k|\langle W_k,t^{\ast}_k \rangle|} \wedge D^{q_k}\|t_k\|^{q_k}_{q_k}.
\end{align}
and finally 
\begin{align}
\nonumber & \E \sup_{t\in S} \sum^M_{k=1} \rbr{Q_k|\langle W_k,t^{\ast}_k \rangle| }\wedge D^{q_k}\|t_k\|^{q_k}_{q_k} \\
\nonumber & \ls \E\sup_{t\in S} \sum^M_{k=1} \rbr{C|\langle W_k,t^{\ast}_k\rangle|} \wedge D^{q_k}\|t_k\|^{q_k}_{q_k}\\
\label{s5-ineq0.6} & +\E \sup_{t\in S} \sum^M_{k=1} \rbr{(Q_k-C)_{+} |\langle W_k,t^{\ast}_k\rangle|}\wedge D^{q_k}\|t_k\|^{q_k}_{q_k}.
\end{align}
We prove that the latter term is small comparable to $p$.   Recall that $|S|\ls |T|\ls  1+e^p$,  we are done if we show that
\begin{align}\label{s5-ineq1}
\left\| \sum^M_{k=1}\rbr{(Q_k-C)_{+} |\langle W_k,t^{\ast}_k\rangle|}\wedge D^{q_k}\|t_k\|^{q_k}_{q_k} \right\|_{p}\ls cp
\end{align}
for $c$ suitably small.  We prove the result in the next theorem.
This is the main estimate in this section.
\bt\label{s5-theo1}
For suitably large $C$, 
\eqref{s5-ineq1} holds with $c$  which can be suitably small.
\et
\begin{dwd}
Clearly,  for $u>\tilde{C}$
\begin{align*}
& \P(Q_k>u+\tilde{C})=\frac{b_k^{\frac{n_k}{p_k}}}{\frac{1}{p_k}\Gamma(\frac{n_k}{p_k})}\int^{\infty}_{\|\tilde{R}_k\|_2(\tilde{C}+u)}x^{n_k-1}e^{-b_kx^{p_k}}dx\\
& =\frac{b_k^{\frac{n_k}{p_k}}}{\frac{1}{p_k}\Gamma(\frac{n_k}{p_k})}\int^{\infty}_{\| \tilde{R}_k\|_2u} 
(\|\tilde{R}_k\|_2 \tilde{C}+x)^{n_k-1}e^{-b_k(\|\tilde{R}_k\|_2\tilde{C}+x)^{p_k}}dx\\
&\ls 2^{n_k-1}\exp(-b_k\| \tilde{R}_k\|^{p_k}_2C^{p_k})\frac{b_k^{\frac{n_k}{p_k}}}{\frac{1}{p_k}\Gamma(\frac{n_k}{p_k})}
\int^{\infty}_{\| \tilde{R}_k\|_2u} 
x^{n_k-1}e^{-b_kx^{p_k}}dx\\
&=2^{n_k-1}\exp(-b_k\| \tilde{R}_k \|_{2}^{p_k}\tilde{C}^{p_k})\P(Q_k>u).
\end{align*}
Therefore,  choosing $\tilde{C}=C/2$,  we have $(Q_k-C)_{+}$ can be replaced by $\delta_k Q_k$, where $\delta_k$ is independent of $Q_k$, moreover
$\delta_k\in \{0,1\}$ and 
\[
\P(\delta_k=1)=2^{n_k-1}\exp(-b_k\|\tilde{R}_k\|_{2}^{p_k}\tilde{C}^{p_k}).
\]
Note that
\[
b_k\|\tilde{R}\|^{p_k}_{2}=\frac{\Gamma\rbr{\frac{n_k+2}{p_k}}^{\frac{p_k}{2}}}{\Gamma\rbr{\frac{n_k}{p_k}}^{\frac{p_k}{2}}}\simeq
\frac{n_k+2}{p_k},
\]
and therefore,  for suitably large $C$, 
\[
\P(\delta_k=1)\ls \exp(-n_k(C/4)^{p_k}).
\]
We denote by $\E_{\delta},\E_{Q,W}$ integration with respect  to variables  $(\delta_k)^M_{k=1}$ and $(Q_k,W_k)^M_{k=1}$ .   We have 
\begin{align}
\nonumber & \|\sum^M_{k=1}\rbr{(Q_k-C)_{+} |\langle W_k,t^{\ast}_k\rangle|}\wedge D^{q_k}\|t_k\|^{q_k}_{q_k}\|_p\ls \|\sum^M_{k=1} \delta_k \rbr{Q_k|\langle  W_k,t^{\ast}_k \rangle|}  \wedge D^{q_k}\|t_k\|^{q_k}_{q_k} \|_p\\
\label{s5-ineq1.5} &=\sbr{\E_{\delta} \sum_{K\subset [M]} \prod_{k\in K} 1_{\delta_k=1}\prod_{l\in K^c}1_{\delta_l=0} \E_{Q,W} \sbr{ \sum_{k\in K} \rbr{Q_k |\langle W_k,t_k\rangle|}
\wedge D^{q_k}\|t_k\|^{q_k}_{q_k}}^p}^{\frac{1}{p}}.
\end{align}
We use Proposition \ref{s3-prop2} - the equation \eqref{s3-ineq13.5} to get
\begin{align*}
 & \rbr{\E_{Q,W} \sbr{ \sum_{k\in K} \rbr{Q_k |\langle W_k,t^{\ast}_k\rangle|}
\wedge D^{q_k}\|t_k\|^{q_k}_{q_k}}^p}^{\frac{1}{p}}\ls \|\sum^M_{k=1}Q_k|\langle W_k,t^{\ast}_k\rangle| \wedge  D^{q_k}\|t_k\|^{q_k}_{q_k} \|_p\\
 & \lesssim \sum^M_{k=1} \rbr{r_k^{\frac{1}{p_k}}\|t^{\ast}_k\|_{q_k}\wedge D^{q_k}\|t_k\|^{q_k}_{q_k}}1_{r_k>|I_k(t)|}
\end{align*}
for some $\sum^M_{k=1}r_k=p$.  Let us denote by $K_0$ the subset of $[M]$ that consists of $k$
such that $r_k^{1/p_k}\|t^{\ast}_k\|_{q_k}\ls Dr_k$.  Clearly,
\[
\sum_{k\in K_0}  r_k^{\frac{1}{p_k}}\|t^{\ast}_k\|_{q_k}\wedge D^{q_k}\|t_k\|^{q_k}_{q_k}\ls D\sum_{k\in K_0}r_k \lesssim Dp.
\]
On the other hand,  if for some $k\in [M]$,  $r_k^{1/p_k}\|t^{\ast}_k\|_{q_k}> Dr_k$,  then obviously    $t_k=t_k^{\ast}$ and
 $r_k\ls D^{q_k}\|t_k\|_{q_k}^{q_k}$.  However,  in the case of small coefficients \eqref{s4-ineq2.1} we have $D^{q_k}\|t_k\|_{q_k}^{q_k}\ls A^{p_k}n_k$
and hence $r_k\ls  A^{p_k}n_k$ for all $k\in K_0^{c}.$.  Once again by Proposition \ref{s3-prop2} - the equation \eqref{s3-ineq13}
\[
\sum_{k\in K^c_0}r_k^{1/p_k}\|t^{\ast}_k\|_{q_k}1_{r_k>|I_k(t)|} \lesssim A\|\sum^M_{k=1}|P_k |\langle W_k,t^{\ast}_k\rangle|\|_p \lesssim A\|X_{t^{\ast}}\|_p \lesssim Ap,
\]
where we have used the simplification from Proposition \ref{s2-prop2} - i.e.  $\|X_t\|_p\ls 2p$ and the Bernoulli comparison which gives $\|X_{t^{\ast}}\|_p\ls \|X_t\|_p$.  In this way we get
\[
 \rbr{\E_{Q,W} \sbr{ \sum_{k\in K} \rbr{Q_k |\langle W_k,t^{\ast}_k\rangle|}
\wedge D^{q_k}\|t_k\|^{q_k}_{q_k}}^p}^{\frac{1}{p}}\lesssim \max\{A,D\}p.
\]
Let us return to our bound \eqref{s5-ineq1.5}.
Note that we should only care for $K\subset M$ such that
\[
\| \sum_{k\in K}\rbr{ Q_k |\langle W_k,t^{\ast}_k\rangle|} \wedge D^{q_k}\|t_k\|_{q_k}^{q_k}\|_{p}\gs pc.
\]
We use now that  $D^{q_k}\|t^{\ast}_k\|^{q_k}_{q_k}\ls A^{p_k}n_k$, 
and hence  $\sum_{k\in K} A^{p_k}n_k \gs pc$.  But then
\[
\E_{\delta} \prod_{k\in K} 1_{\delta_k=1}=\exp(-\sum_{k\in K} (C/4)^{p_k}n_k)\ls \exp(-(C/(4A))cp)\ls \exp(-p/c),
\]
whenever  $C\gs 4A/c^2$.  Finally,  we should observe that due to \eqref{s2-ineq10} we have $|I(t)|\ls \delta' p$. and therefore there are at most $2^{p\delta'}$ sets $K$ for which we have to use the second method.  With respect to  \eqref{s5-ineq1.5} the above bounds imply
\[
\|\sum^M_{k=1}\rbr{(Q_k-C)_{+} |\langle W_k,t^{\ast}_k\rangle|}\wedge D^{q_k}\|t_k\|^{q_k}_{q_k}\|_p \lesssim cp+2^{\delta'}e^{-\frac{1}{c}}\max\{A,D\}p\lesssim cp
\]
if $c$ is suitably small.  It proves the result.
\end{dwd}
Consequently,  by \eqref{s1-rema1}
\[
\E \sup_{t\in S} \sum^M_{k=1} \rbr{(Q_k-C)_{+} |\langle W_k,t^{\ast}_k\rangle|}\wedge D^{q_k}\|t_k\|^{q_k}_{q_k}\ls ecp.
\]
which together with \eqref{s5-ineq0.5} and \eqref{s5-ineq0.6} gives
\[
\E\sup_{t\in S} \sum^M_{k=1} |\langle W_k,t^{\ast}_k\rangle| \wedge D^{q_k}\|t_k\|^{q_k}_{q_k} \gs C^{-1}(8K)^{-1} p,
\] 
for $c$ suitably small,  i.e.  $ec\ls (8K)^{-1}$.  Hence $\E\sup_{t\in S}|\bar{X}|_t\gs (8CK)^{-1}p$.  
As stated in \eqref{s4-ineq4.5} we have solved the case of small coefficients (cf.  \eqref{s4-ineq2.1}).

\section{Large coefficients}

We work in the cube-like setting formulated in Proposition \ref{s2-prop1} accompanied by  the simplification from Proposition \ref{s2-prop2}.
Our goal is to prove the minoration  for large coefficients - see \eqref{s4-ineq2.2}.   Recall now that we can work under the condition 
 \eqref{s4-ineq4}.  Therefore we may assume  that there is a point $t_0\in T$ and a  subset $S\subset T$ such that $|S|\gs e^{p/2}$ and
\[
\sum^M_{k=1} \sum_{i\in I_k(t^{\ast})\triangle I_k(t_{0}^{\ast})} D^{q_k}k_i^{q_k}< p. 
\]
We are going to show that in this case  $\|X_{t^{\dagger}}-X_{s^{\dagger}}\|_p\gtrsim p$ for all $s\neq t$,  $s,t\in S$.
It suffices to prove that $\|X_{t^{\ast}}-X_{t^{\ast}_0}\|_p$ and $\|X_{s^{\ast}}-X_{t^{\ast}_0}\|_p$ are bit smaller than $p$.
More precisely, 
\begin{lema}\label{s6-lema1}
If $D$ is suitably large then  $\|X_{t^{\ast}}-X_{t^{\ast}_0}\|_p\ls p/4$ for any $t\in S$.
\end{lema}
\begin{dwd}
Consider $\sum^M_{k=1}r_k=p$ and
\[
 \|X_{t^{\ast}}-X_{t_0^{\ast}}\|_p\simeq \sum^M_{k=1}\|P_k\|_{r_k}\|\langle W_k,t^{\ast}_k-t^{\ast} _{0k}\rangle \|_{r_k}.
\] 
Not that we have the inequality 
\[
\|P_k\|_{r_k}\\|\langle W_k,t^{\ast}_k-t^{\ast} _{0k}\rangle \|_{r_k}\lesssim
r_k^{\frac{1}{p_k}}\|t^{\ast}_k-t^{\ast} _{0k} \|_{q_k}.
\] 
Using that $a^{1/p_k}b^{1/q_k}\ls a/p_k+ b/q_k$ we get
\[
\|P_k\|_{r_k}\\|\langle W_k,t^{\ast}_k-t^{\ast} _{0k}\rangle \|_{r_k}\ls p_k^{-1}\frac{r_k}{D}+ q_k^{-1}D^{\frac{q_k}{p_k}}\|t^{\ast}_k-t^{\ast}_{0k}\|^{q_k}_{q_k}.
\]
Clearly $q_k/p_k=1/(p_k-1)$.   Now we can benefit from the condition \eqref{s4-ineq4},  namely we have  
\[
\sum^M_{k=1} \sum_{i\in I_k(t^{\ast})\triangle I_k(t_{0}^{\ast})} D^{\frac{1}{p_k-1}}k_i^{q_k}<p/D. 
\]
For $D$ suitably small it proves that $\|X_{t^{\ast}}-X_{t_0^{\ast}}\|_p\ls p/4$.
\end{dwd}
Consequently,  since $\|X_t - X_s\|_p\gs p$ for $s\neq t$,  $s,t\in S$ we obtain
\begin{align}\label{s6-ineq0}
\|X_{t^{\dagger}}-X_{s^{\dagger}}\|_{p}\gs \|X_t-X_{s}\|_p -\|X_{t^{\ast}}-X_{t_0^{\ast}}\|_p-\|X_{s^{\ast}}-X_{t_0^{\ast}}\|_p\gs    p/2.
\end{align}
The last result we need
stems from the fact that $\|X_{t^{\dagger}}\|_p\ls \|X_t\|_p\lesssim p$,  which is due to Bernoulli comparison and Proposition \ref{s2-prop2}.
Namely,  we have
\begin{lema}\label{s6-lema2}
For any $t\in T$,  the following inequality holds
\[
\sum^{M}_{k=1}n_k^{\frac{1}{p_k}}\|t^{\dagger}_k\|_{q_k} \lesssim p.
\]
Moreover,  $\sum^M_{k=1}n_k 1_{t^{\dagger}_k\neq 0}\ls pc$,  where $c$ can be suitably small.
\end{lema}
\begin{dwd}
Obviously $\|X_{t^{\dagger}}\|_p=\|\sum^M_{k=1}P_k\langle W_k,t^{\dagger}_k\rangle\|_p$,  so since $\|X_{t^{\dagger}}\|_p\lesssim p$
there must exist $r_k$,  $1\ls k\ls M$ such that $\sum^M_{k=1}r_k=p$ and
\[
\sum^M_{k=1}\|P_k\|_{r_k}\|\langle W_k,t^{\dagger}_k \rangle \|_{r_k} \lesssim p.
\]
If $r_k\gs n_k$ we can use
\[
\|P_k\|_{r_k}\|\langle W_k,t^{\dagger}_k \rangle \|_{r_k} \gtrsim n_k^{\frac{1}{p_k}}\|t^{\dagger}_k\|_{q_k}.
\]
Our aim is to show that we can use $r_k\gs n_k$ for all $k$.  This possible if
\[
\sum^M_{k=1} n_k1_{t^{\dagger_k\neq 0}}\ls p.
\]
Suppose conversely that  $\sum_{k\in J} n_k 1_{t^{\dagger}_k\neq 0}\simeq p$,  for some $J\subset [M]$ then
the inequality $\|X_{t^{\dagger}}\|_p\lesssim p$ gives
\[
\sum_{k\in J} n_k^{\frac{1}{p_k}}\|t^{\dagger}_k\|_{q_k}\lesssim p.
\]
However,  if $t^{\dagger}_k\neq 0$,  then  $D\|t^{\dagger}_k\|_{q_k}> A^{p_k-1}n_k^{\frac{1}{q_k}}$  by \eqref{s4-ineq2.2}.
This is exactly the moment,  where we need our technical assumption that $p_k$ are not too close to $1$.
Namely,  if $A\gtrsim \max\{ D^{\frac{1}{p_k-1}}\}$ we get  $\sum_{k\in J}\sum^M_{k=1} n_k1_{t^{\dagger}_k\ neq 0}>p$
which is a contradiction.  Consequently,  we can select $r_k\gs n_k$ for any $k$ such that $t^{\dagger}_k\neq 0$.
Thus finally,
\[
\sum^M_{k=1} n_k^{\frac{1}{p_k}}\|t^{\dagger}_k\|_{q_k}\lesssim p.
\]
Moreover,
\[
\sum^M_{k=1} n_k \frac{A^{p_k-1}}{D}\lesssim p.
\]
Since $A$ can be much larger than $\max{D^\frac{1}{p_k-1}}$ it completes the proof.
\end{dwd}
We can start the main proof in the case of large coefficients - cf.  \eqref{s4-ineq2.2}.
We are going to use a similar trick to the approach presented in the 'simplification lemma' - Proposition \ref{s2-prop1}.
The point is that on each $\R^{n_k}=\R^{|J_k|}$ we impose a radial-type distribution.
Namely,  let $\mu_k$ be the probability distribution on $\R^{n_k}$ with the density
\[
\mu_k(dx)=\exp\rbr{-\frac{1}{B}n_k^{\frac{1}{p_k}}\|x\|_{q_k}} \frac{n_k^{\frac{n_k}{p_kq_k}}}{B^{\frac{n_k}{q_k}}\Gamma\rbr{\frac{n_k}{q_k}}|\partial B^{n_k}_{q_k}|}dx.
\]
The fundamental property of $\mu_k$ is that
\[
\mu_k(n^{-\frac{1}{p_k}}BuB^{n_k}_{q_k})=\frac{1}{\Gamma\rbr{\frac{n_k}{q_k}}} \int^u_0 s^{\frac{n_k}{q_k}-1}e^{-s}ds.
\]
The median value of $u_0$ the distribution $\Gamma\rbr{\frac{n_k}{q_k},1}$ is comparable to $\frac{n_k}{q_k}$ and hence 
for  $u= \rho \frac{n_k}{q_k}$,  where $\rho$ is smaller than  $1$   we have 
\[
\frac{1}{\Gamma\rbr{\frac{n_k}{q_k}}}\int^u_0 s^{\frac{n_k}{q_k}-1}e^{-s}ds\gtrsim \rho^{\frac{n_k}{q_k}}\gs  \frac{1}{2} \rbr{\frac{u}{u_0}}^{\frac{n_k}{q_k}}.
\]
Since $u_0\ls n_k/q_k $ it gives
\[
\mu_k(\rho B n_k^{\frac{1}{q_k}}q_k^{-1})\gs  \rho^{\frac{n_k}{q_k}},
\]
for any $\rho$ which is bit smaller that $1$.
Furthermore,  by the construction
\[
\mu_k\rbr{ t^{\dagger}_k+ \rho B n_k^{\frac{1}{q_k}}q_k^{-1} B^{n_k}_{q_k} }\gtrsim  \rho^{\frac{n_k}{q_k}} e^{-B^{-1}n_k^{\frac{1}{p_k}}\|t^{\dagger}_k\|_{q_k}}=\exp\rbr{-\frac{n_k}{q_k}\log\frac{1}{\rho}-B^{-1}n_k^{\frac{1}{p_k}}\|t^{\dagger}_k\|_{q_k} }.
\]
Due to Lemma \ref{s6-lema2} we may find constant $B\gs 1$ in such a way that for any $t\in T$, 
\begin{align}\label{s6-ineq1}
B^{-1}\sum^M_{k=1}n_k^{\frac{1}{p_k}}\|t^{\dagger}_k\|_{q_k}\ls p/8.
\end{align}
Moreover, by the same result we may require that 
\begin{align}\label{s6-ineq2}
 \log\rbr{\frac{1}{\rho}}\cdot  \sum^M_{k=1}n_k 1_{t^{\dagger}_k\neq 0}\ls p/8
\end{align}
for some suitably small $\rho$.
\smallskip

\noindent
Let $\mu$ be a measure defined on $\R^d$ by $\mu=\mu_1\otimes \mu_2 \otimes \ldots \otimes \mu_M$. 
We define also
\[
A_{t}=\rbr{x\in \R^d:\;\; \|t^{\dagger}_k-x_k\|_{q_k}\ls \rho B n_k^{\frac{1}{q_k}}q_k^{-1},\;\;\mbox{or}\;\; \|t^{\dagger}_k\|_{q_k}\ls \rho B n_k^{\frac{1}{q_k}}q_k^{-1},  \;\;\mbox{for all}\;\;k\in [M]}.
\]
Using \eqref{s6-ineq1} and \eqref{s6-ineq2}, we get
\[
\mu(A_t)\gs \exp\rbr{-\sum^M_{k=1}\rbr{\frac{n_k}{q_k}\log\frac{1}{\rho}1_{t^{\dagger}_k\neq 0}+n_k^\frac{1}{p_k}\|t^{\dagger}_k\|_{q_k}} }\gs e^{-p/4}.
\]
However,  there are at least $e^{p/2}$ points in $S$ and hence
\[
\sum_{t\in S} \mu(A_t)\gs e^{p/4}.
\]
We define
\[
A_x=\rbr{t\in T:\;\; \|t^{\dagger}_k-x_k\|_{q_k}\ls \rho B n_k^{\frac{1}{q_k}}q_k^{-1},\;\;\mbox{or}\;\; \|t^{\dagger}_k\|_{q_k}\ls \rho B n_k^{\frac{1}{q_k}}q_k^{-1},  \;\;\mbox{for all}\;\; k\in[M] }
\]
and observe that
\begin{align*}
&\int |A_x| \mu(dx)= \sum_{t\in T} \int 1_{t\in A_x} \mu(dx)=\sum_{t\in T} \int 1_{x\in A_t}\mu(dx)\\
&= \sum_{t\in T} \mu(A_t)\gs e^{p/4} .
\end{align*}
Consequently,  there must exists $x\in \R^d$ such that $|A_x|$ counts at least $e^{p/4}$
points.  Note that we can  require that $\|x_k\|_{q_k}>\rho B n_k^{\frac{1}{q_k}}q_k^{-1}$.  Let us denote the improved set $S$ by $\bar{S}$,  in particular 
$|\bar{S}|>e^{p/4}$.
\smallskip

\noindent The final step is as follows. 
 We define the function 
\[
\varphi_k(t)=\left\{ \begin{array}{lll} x_k\in \R^{n_k} & \mbox{if} & \|t_k\|_{q_k}> \rho B n^{\frac{1}{q_k}}q_k^{-1},\\
 0\in \R^{n_k} & \mbox{if} & \|t_k\|_{q_k}\ls \rho B n_k^{\frac{1}{q_k}}q_k^{-1} \end{array}\right. 
\]
Let $\varphi=(\varphi_k)^M_{k=1}$.  Note that if $\varphi_k(t)=0$ then $t^{\dagger}_k=0$ by  \eqref{s4-ineq2.2}.
Therefore,  we are sure that $\|t^{\dagger}_k-\varphi_k(t)\|_{q_k}\ls \rho B n_k^{\frac{1}{q_k}}q_k^{-1}$.
\begin{lema}\label{s6-lema3}
The following inequality holds
\[
\|X_{t^{\dagger}}-X_{\varphi(t)}\|_p\lesssim cp,
\]
where $c$ is suitably small.
\end{lema}
\begin{dwd}
Indeed it suffices to check for $\sum^M_{k=1}r_k=p$ that
\[
\sum^M_{k=1}r_k^{\frac{1}{p_k}}\|t^{\dagger}_k-\varphi_k(t)\|_{q_k}\ls cp.
\]
We use the inequality $a^{1/p_k}b^{1/q_k}\ls a/p_k+ b/q_k$ and $\|t_k-\varphi_k(t)\|_{q_k}\ls \rho B n_k^{\frac{1}{q_k}}q_k^{-1}$  to get 
\[
r_k^{\frac{1}{p_k}}\|t^{\dagger}_k-\varphi_k(t)\|_{q_k}\ls \frac{1}{p_k}cr_k+\frac{1}{q_k}c^{-\frac{1}{p_k-1}}\rho^{q_k} B^{q_k}q_k^{-q_k} n_k1_{t^{\dagger}_k\neq 0}.
\]
It remains to notice that if only $\rho B q_k^{-1}\ls c$
\[
\sum^M_{k=1}\frac{1}{q_k}c^{-\frac{1}{p_k-1}}\rho^{q_k} B^{q_k} q_k^{-q_k} n_k1_{t^{\dagger}_k\neq 0}\ls c\sum^M_{k=1}n_k 1_{t^{\dagger}_k\neq 0}\ls pc.
\]
However,  $\rho$ can be suitably small,  thus the result follows.
\end{dwd}
Now if $c$ is suitably small we get that
\[
\|X_{\varphi(t)}-X_{\varphi(s)}\|_p \gs p/4.
\]
Indeed this is the  consequence of \eqref{s6-ineq0} and $\|X_{\varphi(t)}-X_{t^{\dagger}}\|_p\ls p/8$ and $\|X_{\varphi(s)}-X_{\varphi(t_0)}\|_p\ls p/8$, namely
\[
\|X_{\varphi(t)}-X_{\varphi(s)}\|_p \gs \|X_{t^{\dagger}}-X_{s^{\dagger}}\|_p- \|X_{\varphi(t)}-X_{t^{\dagger}}\|_p-
\|X_{\varphi(s)}-X_{\varphi(t_0)}\|_p.
\]
Moreover,  by Remark \ref{s1-rema1} and $\bar{S}\subset S\subset T$,  $0\in T$ and $|T|\ls 1+e^p$
\[
\E\sup_{t\in \bar{S}}(X_{\varphi(t)}-X_{t^{\dagger}})\ls ecp.
\]
It proves that
\[
\E\sup_{t\in \bar{S}} X_{\varphi(t)} \ls \E\sup_{t\in \bar{S}}X_{t^{\dagger}}+\E\sup_{t\in \bar{S}}(X_{\varphi(t)}-X_{t^{\dagger}})\ls \E\sup_{t\in \bar{S}}X_{t^{\dagger}}+ecp.
\]
It remains to observe that $X_{\varphi_k(t)}$ is either $0$ or $X_{x_k}$.  But vector  $(X_{x_k})^M_{k=1}$,
is  one unconditional,  log-concave  with independent entries,  so we may use the standard Sudakov minoration -
 the main result of \cite{Lat1}.  Therefore,
\[
\E\sup_{t\in \bar{S}} X_{\varphi(t)} \gs \frac{1}{K} p.
\]
It completes the proof of the  Sudakov minoration  in the case of large coefficients (cf.  \eqref{s4-ineq2.2}).
In this way we have completed the program described in Section 2.

\section{The partition scheme}

Having established the Sudakov minoration one can prove that partition scheme and in
this way establish the characterization of $S_X(T)=\E\sup_{t\in T}X_t$.  We need here additionally that \eqref{s1-ineq3}
holds. 
\smallskip

\noindent
Let us recall the general approach to establish the lower bound on $S_X(T)$. 
We define a family of distances $d_n(s,t)=\|X_t-X_s\|_{2^n}$. Moreover, let
$B_n(t,\va)$ be the ball centred at $t$ with radius $\va$ in $d_n$ distance and in the same way 
$\Delta_n(A)$ as the diameter of $A\subset T$ in $d_n$ distance.
By \eqref{s1-ineq1} we know that $d_{n+1}\ls 2 d_n$,  note also that our condition 
\eqref{s1-ineq3} reads as $(1+\va)d_n\ls d_{n+1}$.   Let us recall that the property works if all $p_k$ are smaller than
some $p_{\infty}<\infty$.
Note that this assumption is not easily removable,  since
for Bernoulli canonical processes the theory which we describe below is  does not work.
\smallskip

\noindent
We follow the generic chaining approach for families of distances described in \cite{TG3}. 
Let $N_n=2^{2^n}$, $n\gs 1$, $N_0=1$.
The natural candidate for the family $F_{n,j}$ is $F_{n,j}=F$,  where 
$F(A)=K\E \sup_{t\in A}X_t$,   $A\subset T$ for suitably large constant $K$.  We have to prove the growth condition,  namely
that for some $r=2^{\kappa-2}\gs 4$ and fixed $n_0\gs 0$ for any given $n\gs n_0$
and $j\in \Z$ if we can find points $t_1,t_2,\ldots t_{N_n}\in B_n(t,2^nr^{-j})$
which are $2^{n+1}r^{-j-1}$ separated in $d_{n+1}$ distance, then for any
sets $H_i\subset B_{n+\kappa}(t_i,2^{n+\kappa}r^{-j-2})$ the following inequality
holds
\begin{align}\label{s7-ineq0}
F\rbr{\bigcup^{N_n}_{i=1}H_i}\gs 2^nr^{-j-1}+\min_{1\ls i\ls N_n}F(H_i).
\end{align}
Note that in our setting
\begin{align}\label{s7-ineq1}
B_n(t,2^nr^{j-2})\subset B_{n+1}(t,2^{n+1}r^{j-2})\subset \ldots\subset 
B_{n+\kappa}(t,2^{n+\kappa}r^{-j-2})
\end{align}
and also
\begin{align}\label{s7-ineq2}
B_{n+\kappa}(t,2^{n+\kappa}r^{-j-2})\subset 
B_n\rbr{t,\frac{2^{n+\kappa}r^{-j-2}}{(1+\va)^{\kappa}}}\subset B_n\rbr{t,\frac{2^{n+2}r^{-j-1}}{(1+\va)^{\kappa}}}.
\end{align}
In particular we get from \eqref{s7-ineq2} that
\begin{align}\label{s7-ineq2.5}
H_i\subset B_{n+\kappa}(t_i,2^{n+\kappa}r^{-j-2})\subset B_n\rbr{ t_i,\frac{2^{n+2}r^{-j-2}}{(1+\va)^{\kappa}}}
\end{align}
and hence we have that $H_i$ are small when compared with the separation level $r^{-j-1}$.
Clearly if $d_{n+1}(t_i,t_j)\gs 2^{n+1}r^{-j-1}$ then
\[
d_n(t_i,t_j)\gs 2^{-1}d_{n+1}(t_i,t_j)\gs 2^nr^{-j-1}.
\]
The last property we need is a special form of concentration,  i.e.
\begin{align}\label{s7-ineq3}
\| (\E \sup_{t\in A}X_t-\sup_{t\in A} X_t)_{+}\|_p\ls L\sup_{t\in A}\|X_t\|_p.
\end{align}
Let us recall that due to the result of \cite{La-Wo} and some straightforward observations fortunately the 
inequality holds in our setting.  We turn to prove the main result of this section.
\begin{prop}\label{s7-theo2}
Suppose that $X_t$ satisfies assumptions on the regularity of distances
as well as the concentration inequality \eqref{s7-ineq3} then
$F_{n,j}=F$, $F(A)=K\E \sup_{t\in A}X_t$ satisfies the growth condition \eqref{s7-ineq0}  for 
some $r$ and $n_0$.
\end{prop}
\begin{dwd}
Let us define $A=\bigcup^{N_n}_{i=1}H_i$ and $H_i\subset B_{n+1}(t_i,2^{n+1}r^{-j-1})$.
We have
\begin{align*}
& F(A)\gs K\E \sup_{1\ls i\ls N_n}X_{t_i}+\sup_{t\in H_i}(X_t-X_{t_i})\\
&\gs K\E \sup_{1\ls i \ls N_n}X_{t_i}+\min_{1\ls i\ls N_n}F(H_i-t_i)\\
&-K(\E \sup_{1\ls i\ls N_n}(\E \sup_{t\in H_i}(X_t-X_{t_i})-\sup_{t\in H_i}(X_t-X_{t_i}))\\
&\gs 2^{n+1}r^{-j-1}+\min_{1\ls i\ls N_n}F(H_i)-2KL\sup_{1\ls i\ls N_n}
\| (\E \sup_{t\in H_i}(X_t-X_{t_i})-\sup_{t\in H_i}(X_t-X_{t_i}))_{+}\|_{2^n}\\
&\gs 
2^{n+1}r^{-j-1}+\min_{1\ls i\ls N_n}F(H_i)-2KL\max_{1\ls i\ls N}\Delta_n(H_i), 
\end{align*} 
where we have used here the Sudakov minoration.
Using  \eqref{s7-ineq2.5},  we get
\[
\Delta_n(H_i)\ls  \frac{2^{n+3}r^{j-1}}{(1+\va)^{\kappa}}.
\]
Therefore choosing $\kappa$ sufficiently large,  we can guarantee that
\[
F(A)\gs 2^nr^{-j-1} 
\]
as required.
\end{dwd}
The basic result of \cite{TG3} is that having the growth condition for families of distances,  it is true that
$S_X(T)=\E\sup_{t\in T}X_t$ is comparable with $\gamma_X(T)$,  where
\[
\gamma_X(T)=\inf_{\ccA}\sup_{t\in T}\sum^{\infty}_{n=0}\Delta_n(A_n(t)).
\]
It completes the proof of  Theorem \ref{diety-character}.  Namely,
for some absolute constant $K$,
\[
K^{-1}\gamma_X(T)\ls \E \sup_{t\in T}X_t\gs K \gamma_X(T).
\]  
In this way it establishes a geometric characterization of 
$S_X(T)$,  assuming  our list of conditions.

\end{document}